# STATISTICS OF EXTREMES BY ORACLE ESTIMATION

By Ion Grama and Vladimir Spokoiny

*University of South Brittany and Weierstrass-Institute*

We use the fitted Pareto law to construct an accompanying approximation of the excess distribution function. A selection rule of the location of the excess distribution function is proposed based on a stagewise lack-of-fit testing procedure. Our main result is an oracle type inequality for the Kullback–Leibler loss.

**1. Background and outline of main results.** Let $X_1, \ldots, X_n$, be i.i.d. observations with continuous d.f. $F$ supported on the interval $[x_0, \infty)$, $x_0 \geq 0$. Assume that d.f. $F$ is "heavy tailed," that is, that $F$ belongs to the domain of attraction of the Fréchet law $\Phi_{1/\gamma}(x) = \exp(-x^{-1/\gamma})$, $x \geq 0$, with parameter $1/\gamma$. By Fisher–Trippet–Gnedenko theorem (see Bingham, Goldie and Teugels [2]) this is equivalent to saying that for any $x \geq 1$,

$$(1.1) \qquad F_t(x) \to P_\gamma(x) \qquad \text{as } t \to \infty,$$

where $F_t(x)$ is the excess d.f. over the threshold $t > x_0$ defined by

$$F_t(x) = 1 - \frac{1 - F(xt)}{1 - F(t)}, \qquad x \geq 1,$$

and $P_\theta(x) = 1 - x^{-1/\theta}$, $x \geq 1$, is the standard Pareto d.f. with parameter $\theta > 0$. Relation (1.1) suggests using $P_\gamma(x)$ with estimated $\gamma$ as an approximation of $F_t(x)$ for a given $x$ and large $t$. However, it can be misleading in cases when the convergence to the limit distribution is too slow. This is easily seen by inspecting the trajectories of the *Hill estimator* (1.2) computed from samples drawn from the log-gamma distribution $F(x)$, see Figure 1. It is sometimes called the Hill horror plot, because of the important discrepancy between the Hill estimator and the estimated parameter $\gamma$, even for very large sample sizes (see Embrechts, Klüppelberg and Mikosch [7] or Resnick

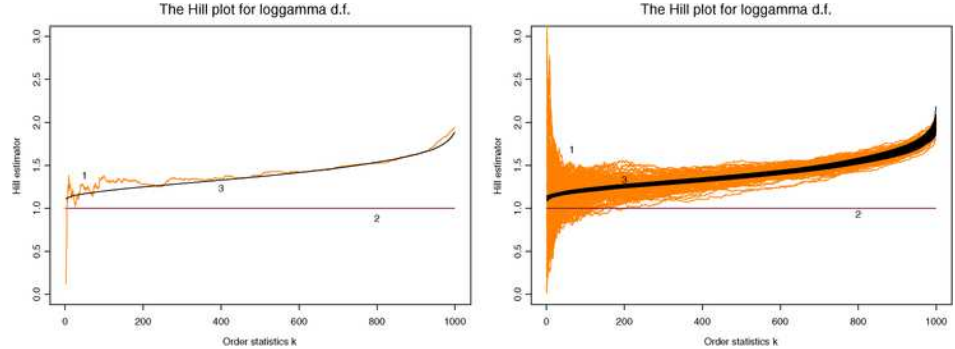

Fig. 1.    1—The Hill estimator $\widehat{h}_{n,k}, k = 1, \ldots, n$, for log-gamma d.f. with rate parameter 1 and shape parameter 2. 2—Index of regular variation $\gamma = 1$ which is expected to be estimated. 3—The fitted Pareto parameter $\theta_t(F)$ computed from the approximation formulas (4.4), (4.8). Left: 1 realization; Right: 100 realizations.

[20]). The explanation lies in the fact that the Hill estimator merely fits a Pareto distribution to the data thereby providing an approximation of the excess d.f. $F_t$ rather than for $\gamma$ itself. Despite these evidences the problem of estimating the excess d.f. $F_t$ regardless of the limit $P_\gamma$ is less studied in the literature.

The goal of the present paper is twofold. First of all, we shall consider the problem of recovering the excess d.f. $F_t$ from the data $X_1, \ldots, X_n$ directly, and second, we shall propose an adaptive procedure of the choice of the location of the tail $t$. Motivated by (1.1), we assume that for large values of $t \geq x_0$ the excess d.f. $F_t$ can be approximated by a Pareto law $P_{\theta_t}$ with some index $\theta_t > 0$ possibly depending on the location $t$ and generally different from $\gamma$. The statistical problem is that of recovering $F_t$ by constructing a family of estimators $\widehat{\theta}_{n,t}, t \geq x_0$, of the parameters $\theta_t, t \geq x_0$, and proposing an adaptive rule for choosing the location threshold $t$.

Some consequences of the main results of the paper are formulated below. Let $X_{n,1} > X_{n,2} > \cdots > X_{n,n}$ be the order statistics pertaining to $X_1, \ldots, X_n$ and $\widehat{h}_{n,k}, k = 1, \ldots, n-1$, be the family of Hill estimators, where

$$(1.2) \qquad \widehat{h}_{n,k} = \frac{1}{k} \sum_{i=1}^{k} \log \frac{X_{n,i}}{X_{n,k+1}},$$

see Hill [14]. Denote $\widehat{n}_t = \sum_{i=1}^{n} 1(X_{n,i} > t)$ and $\widehat{\theta}_{n,t} = \widehat{h}_{n,\widehat{n}_t}$, where $\widehat{\theta}_{n,t} = 0$ if $\widehat{n}_t = 0$.

The discrepancy between two equivalent probability laws $P$ and $Q$ is measured by the Kullback–Leibler divergence $\mathcal{K}(P, Q) = \int \log \frac{dP}{dQ} \, dP$ and by the $\chi^2$-divergence $\chi^2(P, Q) = \int \frac{dP}{dQ} \, dP - 1$. For any $t \geq x_0$ the best approximation of the excess d.f. $F_t$ is defined by looking for the "closest" element in the set of Pareto distributions. Let $\theta_t(F) = \arg\min_{\theta>0} \mathcal{K}(F_t, P_\theta)$ be the minimum Kullback–Leibler divergence Pareto parameter, called in the sequel for short *fitted Pareto index*. Thereafter $\mathbf{P}_F$ denotes the probability measure corresponding to the i.i.d. observations $X_1, \ldots, X_n$ with d.f. $F$.

Instead of (1.1), assume that $F$ admits an accompanying Pareto tail, which means that $\chi^2(F_t, P_{\theta_t(F)}) \to 0$ as $t \to \infty$. This condition is not very restrictive and defines a class of d.f.'s related to those in Hall and Welsh [12] and Drees [6]. Then, according to our Theorem 4.4,

$$\mathcal{K}(F_{\tau_n}, P_{\widehat{\theta}_{n,\tau_n}}) = O_{\mathbf{P}_F}\left(\frac{\log n}{n(1-F(\tau_n))}\right) \qquad \text{as } n \to \infty,$$

for any sequence $\{\tau_n\}$ obeying

$$(1.3) \qquad \chi^2(F_{\tau_n}, P_{\theta_{\tau_n}(F)}) = O\left(\frac{\log n}{n(1-F(\tau_n))}\right) \to 0 \qquad \text{as } n \to \infty.$$

The sequence $\{\tau_n\}$ in the definition of the estimator $\widehat{\theta}_{n,\tau_n}$ being, generally, unknown, we give an automatic selection rule $\widehat{k}_n$ (Section 3) such that the adaptive estimator $\widehat{\theta}_n = \widehat{h}_{n,\widehat{k}_n}$ mimics the nonadaptive estimator $\widehat{\theta}_{n,\tau_n}$, that is,

$$\mathcal{K}(F_{\tau_n}, P_{\widehat{\theta}_n}) = O_{\mathbf{P}_F}\left(\frac{\log n}{n(1-F(\tau_n))}\right) \qquad \text{as } n \to \infty,$$

for any sequence of locations $\{\tau_n\}$ obeying (1.3), see Theorem 4.10. From the results in Hall and Welsh [12] and Drees [6] it follows that the estimators $\widehat{\theta}_{n,\tau_n}$ and $\widehat{\theta}_n$ attain optimal or suboptimal rate in some classes of functions (see Section 5 for details).

Many results on the adaptive choice of the number $k$ of upper statistics involved in the estimation require prior knowledge on the unknown d.f. $F$. A peculiarity of the adaptive procedure proposed in the paper is that it applies to an arbitrary d.f. with Pareto like tail and does not ask additional information on its structure. In particular, $F$ need not even be regularly varying at infinity, that is, it need not satisfy (1.1).

The brief outline of paper: In Section 2 we construct the local likelihood estimators. Section 3 introduces the adaptive procedure for selecting the threshold $t$. Main results of the paper are presented in Section 4. Examples of computing the optimal rates of convergence are given in Section 5. In Sections 7 and 8 we prove exponential type bounds for the likelihood ratio used in the proofs of our main results and necessary auxiliary statements. We shall illustrate the performance of our results on some artificial data sets in Section 6.

**2. Construction of the estimators.** Let $\mathcal{F}$ be the set of all d.f. $F$ having support on the interval $[x_0, \infty)$ with $x_0 \geq 0$, and admitting a strictly positive density $f_F$ w.r.t. Lebesgue measure. For any $t \geq x_0$ define the *excess d.f.* over the threshold $t$ as

$$(2.1) \qquad F_t(x) = 1 - \frac{1 - F(tx)}{1 - F(t)}, \qquad x \geq 1.$$

It is easy to see that

$$F_t(x) = 1 - \exp\left(- \int_t^{tx} \frac{du}{u \alpha_F(u)}\right),$$

where

$$(2.2) \qquad \alpha_F(u) = \frac{1}{u \lambda_F(u)}, \qquad u > x_0,$$

and $\lambda_F(u) = \frac{f_F(u)}{1 - F(u)}, u > x_0$, is the hazard rate. $F_t$ admits the density

$$(2.3) \quad f_{F_t}(x) = \frac{t f_F(tx)}{1 - F(t)} = \frac{1}{x \alpha_F(tx)} \exp\left(- \int_t^{tx} \frac{du}{u \alpha_F(u)}\right), \qquad x \geq 1.$$

Note that, according to the von Mises theorem, if there exists a constant $\alpha > 0$ such that $\alpha_F(x) \to \alpha$ as $x \to \infty$, then $F$ is regularly varying with the index of regular variation $\alpha$, see Beirlant et al. [1].

Recall that given $X_{n,k+1} = t$ the observations $X_{n,1}/t, \ldots, X_{n,k}/t$ are the order statistics of an i.i.d. sequence with common density $f_{F_t}$ (see Reiss [19]). Motivated by this we define the *local* log-*likelihood function*

$$(2.4) \qquad L_{n,t}(F) = \sum_{i: X_i > t} \log f_{F_t}(X_i/t).$$

Let $\mathcal{K}(\theta', \theta) = \mathcal{K}(P_{\theta'}, P_\theta)$ be the Kullback–Leibler divergence between $P_\theta$ and $P_{\theta'}$,

$$(2.5) \qquad \mathcal{K}(\theta', \theta) = \int \log \frac{dP_{\theta'}}{dP_\theta} \, dP_{\theta'} = G\left(\frac{\theta'}{\theta} - 1\right), \qquad \theta', \theta > 0,$$

where $G(x) = x - \log(1 + x)$. We extend this definition by setting $\mathcal{K}(\theta', \theta) = \infty$ if at least one of $\theta' = 0$ or $\theta = 0$ holds. Lemma 8.1 implies

$$\mathcal{K}(\theta_1, \theta_2) \asymp \left(\frac{\theta_1}{\theta_2} - 1\right)^2 \qquad \text{as } \frac{\theta_1}{\theta_2} - 1 \to 0.$$

2.1. *Pareto-type tails.* Let $\mathcal{F}_t$ be the set of functions $F \in \mathcal{F}$ satisfying $\alpha_F(x) = \theta$, for $x \in (t, \infty)$, where $\theta > 0$ and $t \geq x_0$. If $F \in \mathcal{F}_t$, then the d.f. $F_t$ is exactly Pareto $P_\theta$. Maximization of the local log-likelihood (2.4) over $\mathcal{F}_t$ gives the maximum local quasi-likelihood estimator

$$(2.6) \qquad \widehat{\theta}_{n,t} = \frac{1}{\widehat{n}_t} \sum_{i \,:\, X_i \in (t,\infty)} \log \frac{X_i}{t},$$

where $\widehat{n}_t = \sum_{i=1}^{n} 1(X_i > t)$ denotes the number of observations in the interval $(t, \infty)$. Here and in the sequel the indeterminacy $0/0$ arising in the definition of the estimators is understood as $0$, that is, for $t \geq X_{n,1}$ the estimator $\widehat{\theta}_{n,t}$ is defined to be $0$. Although $\widehat{\theta}_{n,t}$ is not exactly the Hill estimator, it is closely related. In fact, if $t = X_{n,k+1}$, where $1 \leq k \leq n-1$, then $\widehat{\theta}_{n,t} = \widehat{\theta}_{n,X_{n,k+1}}$ coincides with the Hill estimator $\widehat{h}_{n,k}$, see (1.2).

Let $L_{n,t}(\theta', \theta) = L_{n,t}(P_{\theta'}) - L_{n,t}(P_\theta)$ be the log of the local likelihood ratio of $P_{\theta'}$ w.r.t. $P_\theta$. By elementary calculation one can see that

$$(2.7) \qquad L_{n,t}(\widehat{\theta}_{n,t}, \theta) = \widehat{n}_t \mathcal{K}(\widehat{\theta}_{n,t}, \theta).$$

2.2. *Pareto change point-type tails.* Let $\mathcal{F}_{t,\tau}$ be the set of functions $F \in \mathcal{F}$ having the change point structure: $\alpha_F(x) = \theta_1$, for $x \in [t, \tau)$, $\alpha_F(x) = \theta_2$, for $x \in [\tau, \infty)$, where $\theta_1, \theta_2 > 0$ and $1 \leq t \leq \tau < \infty$. Of course $\mathcal{F}_t \subset \mathcal{F}_{t,\tau}$. If $F \in \mathcal{F}_{t,\tau}$, then the d.f. $F_t$ coincides with the Pareto change point d.f.

$$P_{\theta_1, \theta_2, \tau/t}(x) = 1 - \exp\left( \int_1^x \frac{du}{\alpha'(u)u} \right),$$

where $\alpha'(x) = \theta_1$, for $x \in [1, \tau/t)$, $\alpha'(x) = \theta_2$, for $x \in [\tau/t, \infty)$. For given $t \leq X_{n,1}$ and $\tau \geq 1$ maximization of the local likelihood (2.4) over $\mathcal{F}_{t,\tau}$ gives the maximum likelihood estimator $(\widehat{\theta}_{n,t,\tau}, \widehat{\theta}_{n,\tau})$, where

$$\widehat{\theta}_{n,t,\tau} = \frac{\widehat{n}_t \widehat{\theta}_{n,t} - \widehat{n}_\tau \widehat{\theta}_{n,\tau}}{\widehat{n}_{t,\tau}}$$

and $\widehat{n}_{t,\tau} = \widehat{n}_t - \widehat{n}_\tau = \sum_{i=1}^{n} 1(t < X_i \leq \tau)$ is the number of observations in the interval $(t, \tau]$. As above, $\widehat{\theta}_{n,t,\tau} = 0$ if $t \geq X_{n,1}$.

Denote by $L_{n,t}(\theta_1, \theta_2, \tau, \theta) = L_{n,t}(P_{\theta_1, \theta_2, \tau/t}) - L_{n,t}(P_\theta)$ the local log-likelihood *ratio* corresponding to Pareto change point model $P_{\theta_1, \theta_2, \tau/t}(x)$ with respect to the Pareto model $P_\theta$. By straightforward calculations it is verified that

$$(2.8) \qquad L_{n,t}(\widehat{\theta}_{n,t,\tau}, \widehat{\theta}_{n,\tau}, \tau, \theta) = \widehat{n}_{t,\tau} \mathcal{K}(\widehat{\theta}_{n,t,\tau}, \theta) + \widehat{n}_\tau \mathcal{K}(\widehat{\theta}_{n,\tau}, \theta).$$

**3. Adaptive selection of the location of the tail.** Several procedures have been proposed in the literature for the choice of the number of upper statistics to be used in the estimation of the index of regular variation. We refer to Beirlant et al. [1], and to the references therein [3, 5, 8, 11, 13, 15]. However one should note that most of these procedures require some prior knowledge on the d.f. $F$.

To illustrate the problem let us recall the main result in Hall and Welsh [12] (see also Drees [6]). Let $F$ be a d.f. with density

$$(3.1) \qquad f_F(x) = d\alpha x^{-(\alpha+1)}(1 + r(x)), \qquad |r(x)| \le Ax^{-\alpha\rho}, x \ge 0,$$

where $|\alpha - \alpha_0| \le \varepsilon, |d - d_0| \le \varepsilon$ and $\alpha_0, d_0, \varepsilon, \rho, A > 0$. It is proved that the optimal rate of convergence that can be achieved for estimating $\alpha = 1/\beta$ is $n^{-\rho/(2\rho+1)}$. This optimal rate is attained for the *Hill estimator* $\widehat{h}_{n,k_n}$ with the choice $k_n \sim n^{2\rho/(2\rho+1)}$ depending on $\rho$. An adaptive estimator can be constructed by estimating $\rho$ and implementing this estimate into the optimal $k_n$. This approach requires us to know in advance the class of distributions $F$, or generally this information is not available in practice. It is also too conservative in the sense that it is oriented to the worst case in the given class but it may happen that particular distributions have nicer properties.

In this paper we will give a selection procedure which is distribution free and attains exactly or nearly optimal rates for each particular law $F$ in contrast to minimax estimation which is oriented to the worst case in a given class of functions. These kinds of results are usually related to the so-called *oracle inequalities* (Donoho and Jonstone [4]).

The selection rule of the location of the tail $\tau$ which we propose is based on the stagewise lack-of-fit testing for the Pareto distribution (see also Grama and Spokoiny [9]). It can be compared with the adaptive procedures for selecting the bandwidth in nonparametric pointwise function estimation, see Lepski [16], Lepski and Spokoiny [17]. Drees and Kaufmann [5] give a variant of the latter adapted to the tail index estimation. A stagewise procedure for testing Pareto d.f. has been proposed and its performance analyzed in Hall and Welsh [13], where it was shown that the choice based on the detection of lack-of-fit point introduces a significant bias. Our procedure differs from these approaches since the point of lack-of-fit serves just as a pilot for the choice of $k$.

3.1. *The lack-of-fit test.* Denote by $[a]$ the integer part of $a$. Assume that the sequence of positive integers $\{K_n\}$ satisfies $K_n \le n$ and $\lim_{n\to\infty} K_n = \infty$. Consider the uniform grid $r_i = r_i(n) = [in/K_n], i = 1, \ldots, K_n$. In particular, if $K_n = n$ we have $r_i = i$, for $i = 1, \ldots, n$. Let $k_0$ be a positive integer much smaller than $n$.

We shall choose the location of the tail of $F$ in the random set $\{X_{n,r_i} : i = k_0, \ldots, K_n\}$ and therefore the problem reduces to the choice of the natural number $r_i$. We shall proceed by local change-point detection, which consists in consecutive testing for the hypothesis $H^0_{n,r_m}$ that conditionally on $X_{n,r_m+1} = s$ the observations $X_{n,1}/s, \ldots, X_{n,r_m}/s$ are the order statistics of an i.i.d. sample with a Pareto d.f. $P_\theta$ against the alternative $H^1_{n,r_m}$ that conditionally on $X_{n,r_m+1} = s$ the observations $X_{n,1}/s, \ldots, X_{n,r_m}/s$ are the order statistics of a i.i.d. sample with a Pareto change-point d.f. $P_{\theta_1, \theta_2, \tau/s}$, for all $m = r_{k_0}, \ldots, r_{K_n}$.

For testing $H^0_{n,r_m}$ against $H^1_{n,r_m}$ we shall make use of the likelihood ratio statistic $T_n(t, \tau)$ which is defined by

$$(3.2) \quad T_n(t, \tau) = \sup_{F \in \mathcal{F}_{t,\tau}} L_{n,t}(F) - \sup_{F \in \mathcal{F}_t} L_{n,t}(F) = L_{n,t}(\widehat{\theta}_{n,t,\tau}, \widehat{\theta}_{n,\tau}, \tau, \widehat{\theta}_{n,t}),$$

for $x_0 \leq t \leq \tau$. Taking into account (2.8) one gets

$$(3.3) \quad T_n(t, \tau) = T_n^{(1)}(t, \tau) + T_n^{(2)}(t, \tau), \qquad t < \tau,$$

where

$$T_n^{(1)}(t, \tau) = \widehat{n}_{t,\tau} \mathcal{K}(\widehat{\theta}_{n,t,\tau}, \widehat{\theta}_{n,t}), \qquad T_n^{(2)}(t, \tau) = \widehat{n}_\tau \mathcal{K}(\widehat{\theta}_{n,\tau}, \widehat{\theta}_{n,t}).$$

For each $m$ and $k \leq m$ consider the test statistics

$$(3.4) \quad T_{n,m} = \max_{\rho m \leq k \leq (1-\delta)m} T_{n,m,k}, \qquad T_{n,m,k} = T_{n,m,k}^{(1)} + T_{n,m,k}^{(2)},$$

where

$$T_{n,m,k}^{(i)} = T_n^{(i)}(X_{n,m}, X_{n,k}), \qquad i = 1, 2,$$

and $\rho$ and $\delta$ are constants satisfying $0 < \rho, \delta \leq \frac{1}{3}$. We shall suppose that $\delta$ is so large that $(1-\delta)r_i \leq r_{i-1}$, for all $i = k_0, \ldots, K_n$. Actually this condition is satisfied for any given $\delta > 0$ when $n$ becomes sufficiently large. We shall also assume that $\rho r_{k_0} \geq r_1$.

The hypothesis $H^0_{n,r_m}$ will be rejected if $T_{n,r_m} > \mathfrak{z}_n$, for some critical value $\mathfrak{z}_n = \mu \log n$, where $\mu$ is a positive constant.

3.2. *The adaptive procedure.* At this stage the required parameters are the number of the points on the grid $K_n$, the starting point $k_0$, two numbers $\rho$ and $\delta$ which determine the size of the testing window and the critical value $\mathfrak{z}_n$.

The procedure of the adaptive choice of the value $\widehat{k}_n$ reads as follows:

**Initialize** Set $i = k_0$.

**Step 1** Compute the test statistic $T_{n,r_i}$ by (3.4).

**Step 2** If $i \leq K_n$ and $T_{n,r_i} \leq \mathfrak{z}_n$, increase $i$ by 1 and repeat the procedure from Step 1. If $i \leq K_n$ and $T_{n,r_i} > \mathfrak{z}_n$, define

$$(3.5) \qquad \widehat{k}_n = \arg \max_{\rho r_i \leq k \leq (1-\delta)r_i} T^{(2)}_{n,r_i,k}$$

and exit the procedure. If $i > K_n$ we define $\widehat{k}_n = n$ and exit the procedure.

The described procedure is equivalent to defining the adaptive value

$$(3.6) \qquad \widehat{k}_n = \arg \max_{\rho \widehat{m}_n \leq k \leq (1-\delta)\widehat{m}_n} T^{(2)}_{n,\widehat{m}_n,k},$$

where

$$(3.7) \qquad \widehat{m}_n = \min\{r_i : T_{n,r_i} > \mathfrak{z}_n, \ i = k_0, \ldots, K_n\},$$

with the convention $\min \varnothing = r_{K_n}$. The adaptive location of the tail $\tau$ is then defined by $\widehat{\tau}_n = X_{n,\widehat{k}_n}$ and the adaptive estimator is set to

$$\widehat{\theta}_n = \widehat{h}_{n,\widehat{k}_n} \equiv \widehat{\theta}_{n,\widehat{\tau}_n}.$$

REMARK 3.1. In the case of Pareto observations the test statistics (3.3) and (3.4) do not depend on the parameter of the Pareto law. This suggests to compute the critical values $\mathfrak{z}_n$ by Monte Carlo simulations from the homogeneous model with i.i.d. standard Pareto observations. Our simulations show that the proposed adaptive procedure is sensitive to some extent to $\rho$, while being less sensitive to $\delta$, $k_0$ and $K_n$. The choice of these parameters is discussed in Section 6. The reason of introducing the parameter $K_n$ is to speed up numerical execution of the adaptive choice. In order to simplify the formulations and the proofs of the results, in the sequel we shall consider only the case $K_n = n$, which means that $\tau$ will be chosen among all order statistics $X_{n,1}, \ldots, X_{n,n}$.

**4. Main results.** Recall that $\widehat{n}_t$ is the number of observations in the interval $(t, \infty)$. Let $n_t = n(1 - F(t))$ be the expected number of observations in the same interval. Note that $\widehat{\theta}_{n,t} = \widehat{h}_{n,\widehat{n}_t}$, $t \geq x_0$, by (1.2) and (2.6).

Thereafter $\mathbf{P}_F$ and $\mathbf{E}_F$ denote the probability and the expectation pertaining to the i.i.d. observations $X_1, \ldots, X_n$ with common d.f. $F$. For any equivalent probability measures $P$ and $Q$ we denoted by $\mathcal{K}(P,Q) = \mathbf{E}_P \log \frac{dP}{dQ}$ the Kullback–Leibler divergence and by $\chi^2(P,D) = \int \frac{dP}{dQ} dP - 1$ the $\chi^2$-divergence. A simple application of Jensen's inequality shows that $0 \leq \mathcal{K}(P,Q) \leq \log(1 + \chi^2(P,Q))$.

We shall measure the discrepancy between two possible values $\theta_1 > 0$ and $\theta_2 > 0$ of the Pareto index in terms of the Kullback–Leibler divergence $\mathcal{K}(\theta_1, \theta_2)$ between two Pareto measures, see (2.5).

4.1. *Rates of convergence of nonadaptive estimators.* We say that the d.f. $F$ admits an accompanying Pareto tail with tail index function $\theta_t$, $t \geq x_0$, if for any $t \geq x_0$ there exists an index $\theta_t > 0$ such that $\theta_t$ is a continuous function of $t$ and

$$\lim_{t \to \infty} \chi^2(F_t, P_{\theta_t}) = 0. \tag{4.1}$$

This definition can be viewed as an extension of the regular variation condition (1.1). Instead of requiring the existence of the limit $P_\gamma$ it stipulates that $F_t$ admits an accompanying Pareto law $P_{\theta_t}$ with a parameter $\theta_t > 0$ possibly changing with $t$. The class of d.f. satisfying (4.1) is very large. For instance the d.f.'s satisfying the Hall condition (3.1), log-gamma d.f. and Pareto d.f. with logarithmic-type perturbations are of this type. We refer to Section 5, where $\theta_t$ is explicitly computed for these examples. The class of distributions defined by (4.1) includes d.f.'s which are not regularly varying. Examples are normal and exponential d.f.'s with some $\theta_t \to 0$ as $t \to \infty$.

It is easy to see that if the d.f. $F$ admits an accompanying Pareto tail with tail index function $\theta_t$, $t \geq x_0$, then there exists a sequence $\{\tau_n\}$ such that

$$\chi^2(F_{\tau_n}, P_{\theta_{\tau_n}}) = O\left(\frac{\log n}{n(1 - F(\tau_n))}\right) \to 0 \qquad \text{as } n \to \infty. \tag{4.2}$$

For the sake of brevity, a sequence of locations $\{\tau_n\}$ satisfying (4.2) is said to be *admissible*.

THEOREM 4.1. *Assume that the d.f. $F$ admits an accompanying Pareto tail with tail index function $\theta_t$, $t \geq x_0$. Then, for any admissible sequence of locations $\{\tau_n\}$,*

$$\mathcal{K}(\widehat{\theta}_{n,\tau_n}, \theta_{\tau_n}) = O_{\mathbf{P}_F}\left(\frac{\log n}{n(1 - F(\tau_n))}\right) \qquad \text{as } n \to \infty.$$

Here and in the sequel the constant in $O_{\mathbf{P}_F}$ depends only on the constant in $O$ in (4.2). This theorem is an immediate consequence of the more general Theorem 4.5 formulated below.

COROLLARY 4.2. *Assume that $F$ admits an accompanying Pareto tail with a constant tail index function $\theta_t = \gamma$, $t \geq x_0$. Then by Theorem 4.1,*

$$\mathcal{K}(\widehat{\theta}_{n,\tau_n}, \gamma) = O_{\mathbf{P}_F}\left(\frac{\log n}{n(1 - F(\tau_n))}\right) \qquad \text{as } n \to \infty, \tag{4.3}$$

*with $\tau_n$ satisfying (4.2).*

For any $t \geq x_0$ we "project" $F_t$ on the set $\mathcal{P}$ by choosing the closest element to $F_t$ in the set of Pareto d.f.'s $\mathcal{P} = \{P_\theta : \theta > 0\}$, say $P_{\theta_t(F)}$, where

$$\theta_t(F) = \arg\min_{\theta > 0} \mathcal{K}(F_t, P_\theta).$$

The parameter $\theta_t(F)$ will be called in the sequel *fitted Pareto index*. It can be easily computed and has the following explicit expression (see Figure 2 for a graphical representation):

$$(4.4) \qquad \theta_t(F) = \int_1^\infty \log x F_t(dx) = \int_t^\infty \log \frac{x}{t} \frac{F(dx)}{1 - F(t)}, \qquad t \geq x_0.$$

In cases when (4.1) holds and $F$ is regularly varying at $\infty$ with index of regular variation $\gamma$, it is easy to verify that $\theta_t(F) \to \gamma$ as $t \to \infty$.

COROLLARY 4.3. *Assume that $F$ admits an accompanying Pareto tail with tail index function $\theta_t = \theta_t(F)$, $t \geq x_0$. Then according to Theorem 4.1*

$$(4.5) \qquad \mathcal{K}(\widehat{\theta}_{n,\tau_n}, \theta_{\tau_n}(F)) = O_{\mathbf{P}_F}\left(\frac{\log n}{n(1 - F(\tau_n))}\right) \qquad \text{as } n \to \infty,$$

*where $\tau_n$ satisfies (4.2).*

Corollaries 4.2 and 4.3 can be compared with the consistency results for the Hill estimator established by Mason [18] (see also Hall [10]). Recall the main result of [18]. If $F$ is regularly varying with index of regular variation $\gamma$ and $k_n$ satisfies $k_n \to \infty$ and $k_n/n \to 0$ as $n \to \infty$, then

$$(4.6) \qquad \widehat{h}_{n,k_n} \overset{\mathbf{P}_F}{\to} \gamma \qquad \text{as } n \to \infty.$$

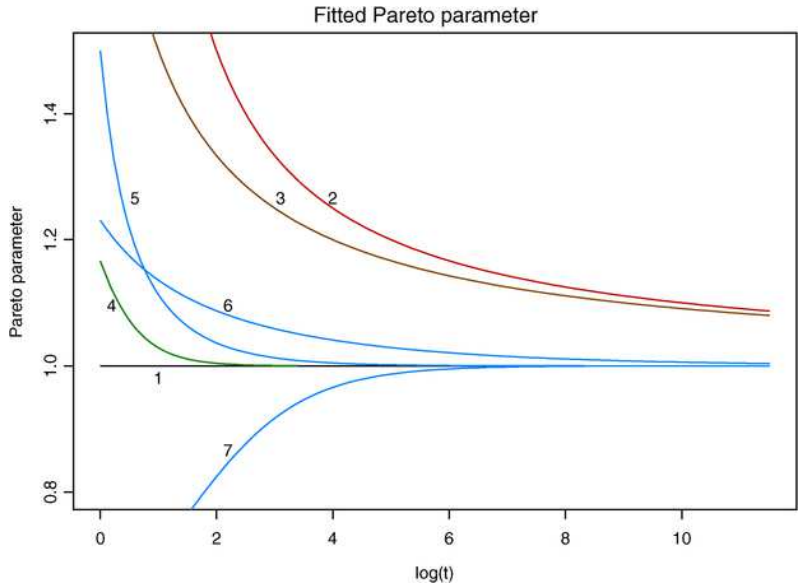

FIG. 2. *Fitted Pareto index $\theta_t(F)$: 1—Pareto d.f.; 2—log-perturbed Pareto d.f. (5.3); 3—log-gamma d.f.; 4—Cauchy d.f.; 5, 6, 7—Hall model (5.4).*

Our Corollary 4.3 improves upon this result by stating that if $F$ admits an accompanying Pareto tail with tail index function $\theta_t = \theta_t(F)$, $t \geq x_0$, then for any $\tau_n$ satisfying (4.2),

$$(4.7) \qquad \widehat{h}_{n,\widehat{n}_{\tau_n}} - \theta_{\tau_n}(F) = \widehat{\theta}_{n,\tau_n} - \theta_{\tau_n}(F) \overset{\mathbf{P}_F}{\to} 0 \qquad \text{as } n \to \infty.$$

A comparison of the precision of the approximations (4.6) and (4.7) is given in Figure 1, where the realizations of the estimator $\widehat{h}_{n,k}$ are plotted as processes in $k$ along with the fitted Pareto index $\theta_t(F)$, for $t = X_{n,1}, \ldots, X_{n,n}$. The underlying d.f. $F(x)$ is the log-gamma one. From these graphs it is seen that for finite sample sizes the Hill estimator $\widehat{h}_{n,k}$ provides a satisfactory approximation of the quantity $\theta_{X_{n,k}}(F)$ while staying far away from the solid straight line corresponding to the parameter of regular variation $\gamma = 1$, except the cases when the fitted Pareto index itself is close to $\gamma$. These conclusions are confirmed also by simulation results reported in Figure 3.

Note that the fitted Pareto index $\theta_t(F)$ coincides with the mean value of the function $\alpha_F$ [see (2.2)] on the interval $[t, \infty)$ w.r.t. $F_t$:

$$\theta_t(F) = \int_1^\infty \alpha_F(tx) F_t(dx) = \int_t^\infty \alpha_F(x) \frac{F(dx)}{1 - F(t)}.$$

For numerical computations of the value $\theta_t(F)$ one can use the following approximation formula:

$$(4.8) \qquad \theta_{X_{n,k}}(F) \approx \frac{1}{k} \sum_{i=1}^k \alpha_F(X_{n,i}).$$

Now we shall present an application of the bound (4.5) to the estimation of the excess d.f. $F_{\tau_n}$.

THEOREM 4.4. *Assume that the d.f. $F$ admits an accompanying Pareto tail with tail index function $\theta_t = \theta_t(F)$, $t \geq x_0$. Then, for any admissible sequence of locations $\{\tau_n\}$,*

$$\mathcal{K}(F_{\tau_n}, P_{\widehat{\theta}_{n,\tau_n}}) = O_{\mathbf{P}_F}\left(\frac{\log n}{n(1 - F(\tau_n))}\right) \qquad \text{as } n \to \infty.$$

PROOF. For any $\theta > 0$ and any $s > x_0$,

$$(4.9) \qquad \mathcal{K}(F_s, P_\theta) = \mathcal{K}(F_s, P_{\theta_s(F)}) + \mathcal{K}(\theta_s(F), \theta).$$

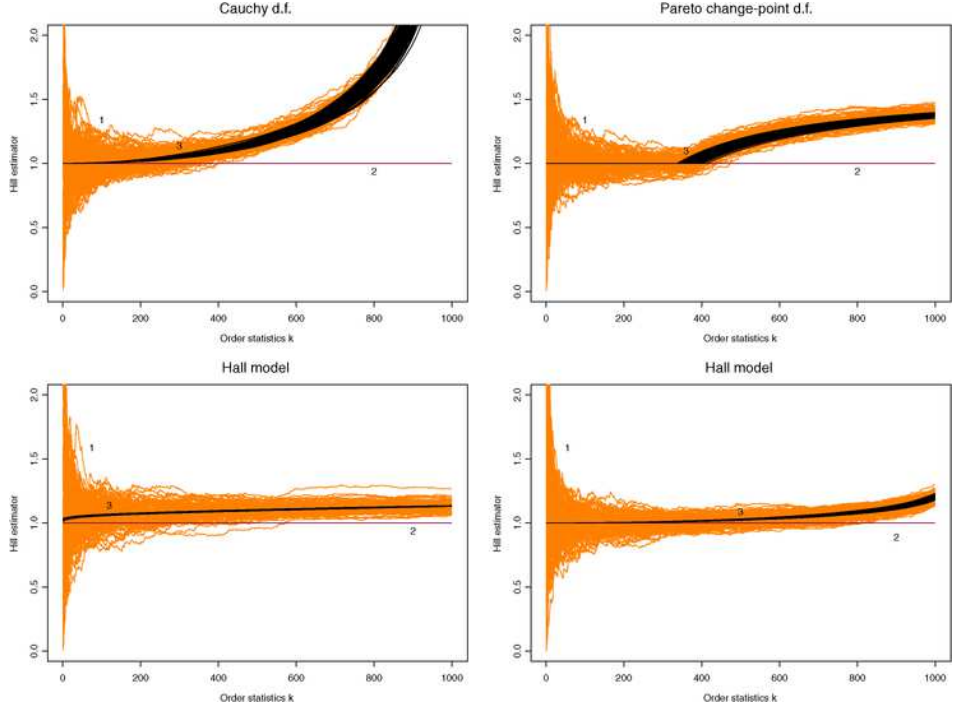

Fig. 3. *1—100 realizations of the Hill estimator* $\widehat{h}_{n,k}, k = 1, \ldots, n,$ *for Cauchy d.f. (top left), Pareto change-point d.f. (top right), Hall model* (5.4) *(bottom left* $\alpha = 1, \beta = 3, c = 1.8$; *bottom right* $\alpha = 1, \beta = 1.2, c = 1.8$); *2—Index of regular variation* $\gamma = 1$ *which is expected to be estimated. 3—The fitted Pareto index* $\theta_1(F)$ *computed from* (4.4), (4.8).

The identity (4.9) follows immediately from the decomposition

$$\mathcal{K}(F_s, P_\theta) = \mathcal{K}(F_s, P_{\theta_s(F)}) + \int_1^\infty \log \frac{dP_{\theta_s(F)}}{dP_\theta} \, dF_s$$

and from

$$\int_1^\infty \log \frac{dP_{\theta_s(F)}}{dP_\theta} \, dF_s = \mathcal{K}(\theta_s(F), \theta).$$

Using (4.9), one gets

$$\mathcal{K}(F_{\tau_n}, P_{\widehat{\theta}_{n,\tau_n}}) = \mathcal{K}(F_{\tau_n}, P_{\theta_{\tau_n}(F)}) + \mathcal{K}(\theta_{\tau_n}(F), \widehat{\theta}_{n,\tau_n}).$$

Since by Lemma 8.1, $\mathcal{K}(\theta_1, \theta_2) \leq \frac{9}{4}\mathcal{K}(\theta_2, \theta_1)$, the assertion follows from the convergence result (4.5) and from the inequality

$$\mathcal{K}(F_{\tau_n}, P_{\theta_{\tau_n}(F)}) \leq \log(1 + \chi^2(F_{\tau_n}, P_{\theta_{\tau_n}(F)})) = O\left(\frac{\log n}{n(1 - F(\tau_n))}\right)$$

as $n \to \infty$. $\square$

The previous results are based on the following more general bound which is a simple application of an exponential bound for the maximum of the likelihood ratio.

THEOREM 4.5. *Assume that $\{\tau_n\}$ is a sequence such that $\tau_n \geq x_0$ and $\lim_{n \to \infty} n(1 - F(\tau_n)) = \infty$. Then for any sequence $\{\theta_n\}$ of positive numbers it holds*

$$\mathcal{K}(\widehat{\theta}_{n,\tau_n}, \theta_n) = O_{\mathbf{P}_F}\left( \frac{\log n}{n(1 - F(\tau_n))} + \chi^2(F_{\tau_n}, P_{\theta_n}) \right) \qquad as \ n \to \infty,$$

*with an absolute constant in $O_{\mathbf{P}_F}$.*

PROOF. Letting $t = s = \tau_n$, $\theta = \theta_n$, $y = 4 \log n + n_{\tau_n} \chi^2(F_{\tau_n}, P_{\theta_n})$, by the first inequality of Proposition 7.3 one gets

$$\mathcal{K}(\widehat{\theta}_{n,\tau_n}, \theta_n) = O_{\mathbf{P}_F}\left( \frac{\log n}{\widehat{n}_{\tau_n}} + \frac{n_{\tau_n}}{\widehat{n}_{\tau_n}} \chi^2(F_{\tau_n}, P_{\theta_n}) \right) \qquad as \ n \to \infty.$$

To finish the proof we use the fact that by Lemma 8.3 it holds $\widehat{n}_{\tau_n} \overset{\mathbf{P}_F}{\asymp} n_{\tau_n}$ as $n \to \infty$, whenever $\lim_{n \to \infty} n_{\tau_n} = \infty$. $\square$

The rate of convergence $\frac{\log n}{n(1 - F(\tau_n))}$ involved in the previous theorems depends on the unknown d.f. $F$ and on the unknown location $\tau_n$. The best possible rate of convergence for a given $F$ is obtained by choosing $\tau_n$ from the balance equation (4.2). Explicit calculation of the resulting rates of convergence for some d.f.'s $F$ are given in Section 5.

4.2. *Stability property of the test statistic.* In the sequel it is assumed that $\mathfrak{z}_n = \mu \log n$, where $\mu > 0$ is a constant. We say that the location $t$ is *accepted* by the testing procedure if $X_{n,r} \geq t$ implies $T_{n,r} \leq \mathfrak{z}_n$. Set $\{t \text{ is accepted}\} \equiv \Omega_{n,t} = \bigcap_{X_{n,r} \geq t}\{T_{n,r} \leq \mathfrak{z}_n\}$.

THEOREM 4.6. *Assume that the d.f. $F$ admits an accompanying Pareto tail with tail index function $\theta_t$, $t \geq x_0$, and $\{\tau_n\}$ is an admissible sequence of locations. Then there exists a finite positive constant $\mu$ such that*

$$\mathbf{P}_F(\tau_n \text{ is accepted}) \equiv \mathbf{P}_F(\Omega_{n,\tau_n}) \to 1 \qquad as \ n \to \infty.$$

PROOF. First note that by Proposition 7.5

$$\mathbf{P}_F\left( \sup_{\tau_n \leq s \leq \tau} T_n(t, \tau) > z \right) \leq 2n^7 \exp(-y/2) + \frac{1}{n} \leq \frac{3}{n},$$

where $y = 16 \log n$ and $z = 2y + 2n_{\tau_n} \chi^2(F_{\tau_n}, P_{\theta_{\tau_n}})$. On the other hand, by (4.2) $z \leq \mathfrak{z}_n$, for some constant $\mu$ and $n$ sufficiently large. Consequently $\Omega_{n,\tau_n}^c \subseteq \{\sup_{\tau_n \leq t \leq \tau} T_n(t, \tau) > z\}$. This implies $\lim_{n \to \infty} \mathbf{P}_F(\Omega_{n,\tau_n}^c) = 0$.  □

REMARK 4.7.    From the preceding proof it can be easily seen that the constant $\mu$ in the definition of the critical value $\mathfrak{z}_n$ depends only on the constant involved in the definition of $O$ in (4.2), say $\lambda$. A simple tracking of constants shows that a crude upper bound for $\mu$ is $32 + 2\lambda e^\lambda$.

4.3. *Rates of convergence of the adaptive estimator.*  First we compare the performance of the adaptive estimator $\widehat{\theta}_n$ with that of the nonadaptive estimator $\widehat{\theta}_{n,\tau_n}$.

THEOREM 4.8.    *Assume that the d.f. $F$ admits an accompanying Pareto tail with tail index function $\theta_t$, $t \geq x_0$, and $\{\tau_n\}$ is an admissible sequence of locations. Then there exists a constant $\mu > 0$ such that*

$$\mathcal{K}(\widehat{\theta}_n, \widehat{\theta}_{n,\tau_n}) = O_{\mathbf{P}_F}\left(\frac{\log n}{n(1 - F(\tau_n))}\right) \qquad as\ n \to \infty.$$

PROOF.    Let $\Omega_{n,\tau_n}^* = \Omega_{n,\tau_n} \cap \{T_{n,k_0} \leq \mathfrak{z}_n\}$. Since by Theorem 4.6 $\mathbf{P}_F(\Omega_{n,\tau_n}) \to 1$ as $n \to \infty$ and by Lemma 8.4 $\mathbf{P}_F(X_{n,k_0} \geq \tau_n) \to 1$ as $n \to \infty$, it holds

$$(4.10) \quad \mathbf{P}_F(\Omega_{n,\tau_n}^*) \geq \mathbf{P}_F(\Omega_{n,\tau_n} \cap \{X_{n,k_0} \geq \tau_n\}) \to 1 \qquad as\ n \to \infty.$$

Denote $\widetilde{m}_n = \widehat{m}_n - 1$. By the definition of $\widehat{m}_n$ on the set $\Omega_{n,\tau_n}^*$ it holds $\widetilde{m}_n \geq \widetilde{n}_{\tau_n}$ (see Section 3.2). We split the further proof into two parts.

First we shall compare $\widehat{h}_{n,\widetilde{m}_n}$ and $\widehat{h}_{n,\widehat{n}_{\tau_n}}$. To this end define the sequence of natural numbers $m_i$, $i = 0, 1, \ldots, i^*$, such that $m_0 = \widetilde{m}_n$ and $m_i$ is the smallest natural number exceeding $m_{i-1}/2$ for $i = 1, 2, \ldots, i^*$, where $i^*$ such that $\rho m_{i^*} \leq \widehat{n}_{\tau_n} \leq (1 - \delta) m_{i^*}$. Let $m_{i^*+1} = \widehat{n}_{\tau_n}$. Since, on the set $\Omega_{n,\tau_n}^*$,

$$(4.11) \quad T_{n,k} \leq \mathfrak{z}_n \equiv \mu \log n \qquad for\ k \in \mathcal{R}_n, k \leq \widetilde{m}_n,$$

by (3.3), with $s = X_{n,m_{i-1}} \leq \tau = X_{n,m_i}$, one gets

$$m_i \mathcal{K}(\widehat{h}_{n,m_{i-1}}, \widehat{h}_{n,m_i}) \leq T_{n,m_{i-1},m_i} \leq \mu \log n, \qquad i = 1, \ldots, i^* + 1,$$

which in turn implies

$$\sum_{i=1}^{i^*} \sqrt{\mathcal{K}(\widehat{h}_{n,m_{i-1}}, \widehat{h}_{n,m_i})} \leq \mu^{1/2} \log^{1/2} n \sum_{i=1}^{i^*} m_i^{-1/2}.$$

Taking into account that $m_i \geq m_{i-1}/2$, for $i = 1, \ldots, i^*$, we obtain

$$\sum_{i=1}^{i^*} m_i^{-1/2} \leq m_{i^*}^{-1/2} \sum_{i=1}^{i^*} 2^{-(i^*-i)/2} \leq 3.5 m_{i^*}^{-1/2}.$$

Since $\widehat{n}_{\tau_n} \le m_{i^*}$, by Lemma 8.2, on the set $\Omega^*_{n,\tau_n}$ it holds

$$(4.12) \qquad \begin{aligned} \sqrt{\mathcal{K}(\widehat{h}_{n,\widetilde{m}_n}, \widehat{h}_{n,\widehat{n}_{\tau_n}})} &\le \frac{3}{2} \sum_{i=1}^{i^*+1} \sqrt{\mathcal{K}(\widehat{h}_{n,m_{i-1}}, \widehat{h}_{n,m_i})} \\ &\le \frac{3 \cdot 4.5}{2} \mu^{1/2} \frac{\log^{1/2} n}{\widetilde{n}_{\tau_n}^{1/2}}. \end{aligned}$$

Now we shall compare $\widehat{h}_{n,\widetilde{m}_n}$ and $\widehat{h}_{n,\widehat{k}_n}$. Recall that by the definition $\widehat{k}_n$ is a natural number satisfying $\rho \widehat{m}_n \le \widehat{k}_n \le (1-\delta) \widehat{m}_n \le \widetilde{m}_n$ (see Section 3.2). Then, on the set $\Omega^*_{n,\tau_n}$, (4.11) implies $T_{n,\widetilde{m}_n,\widehat{k}_n} \le \mathfrak{z}_n$. Since on the same set it holds $\widetilde{m}_n \ge \widehat{n}_{\tau_n}$, we get

$$(4.13) \qquad \sqrt{\mathcal{K}(\widehat{h}_{n,\widehat{k}_n}, \widehat{h}_{n,\widetilde{m}_n})} \le \frac{\mu^{1/2} \log^{1/2} n}{\widetilde{m}_n^{1/2}} \le \mu^{1/2} \frac{\log^{1/2} n}{\widehat{n}_{\tau_n}^{1/2}}.$$

Summing (4.12) and (4.13), by Lemma 8.2 it follows that on the set $\Omega^*_{n,\tau_n}$,

$$(4.14) \qquad \sqrt{\mathcal{K}(\widehat{h}_{n,\widehat{k}_n}, \widehat{h}_{n,\widehat{n}_{\tau_n}})} \le (c\mu)^{1/2} \frac{\log^{1/2} n}{\widehat{n}_{\tau_n}^{1/2}},$$

where $c$ is an absolute constant. Taking into account (4.10),

$$\mathbf{P}_F\left(\mathcal{K}(\widehat{\theta}_n, \widehat{\theta}_{n,\tau_n}) \le c\mu \frac{\log n}{\widehat{n}_{\tau_n}}\right) \to 1 \qquad \text{as } n \to \infty.$$

To get the requested assertion it suffices to replace the random rate of convergence $\widehat{n}_{\tau_n}$ with the deterministic rate $n_{\tau_n} = n(1 - F(\tau_n))$ by Lemma 8.3. $\square$

Combining Theorem 4.8 with Theorem 4.1 one gets the following assertion:

THEOREM 4.9. *Assume that the d.f. $F$ admits an accompanying Pareto tail with tail index function $\theta_t$, $t \ge x_0$, and $\{\tau_n\}$ is an admissible sequence of locations. Then there exists a constant $\mu > 0$ such that*

$$\mathcal{K}(\widehat{\theta}_n, \theta_{\tau_n}) = O_{\mathbf{P}_F}\left(\frac{\log n}{n(1 - F(\tau_n))}\right) \qquad \text{as } n \to \infty.$$

In particular if condition (4.1) is fulfilled with $\theta_t = \theta_t(F)$ one gets

$$(4.15) \qquad \mathcal{K}(\widehat{\theta}_n, \theta_{\tau_n}(F)) = O_{\mathbf{P}_F}\left(\frac{\log n}{n(1 - F(\tau_n))}\right) \qquad \text{as } n \to \infty.$$

Another case of interest is when $F$ is regularly varying with index of regular variation $\gamma > 0$. Assume that condition (4.1) is satisfied with $\theta_t = \gamma$, $t \geq x_0$. Then

$$(4.16) \qquad \mathcal{K}(\widehat{\theta}_n, \gamma) = O_{\mathbf{P}_F}\left(\frac{\log n}{n(1 - F(\tau_n))}\right) \qquad \text{as } n \to \infty.$$

Now we are in position to formulate the result concerning the approximation of the excess d.f. $F_{\tau_n}$.

THEOREM 4.10. *Assume that the d.f. $F$ admits an accompanying Pareto tail with tail index function $\theta_t = \theta_t(F)$, $t \geq x_0$, and $\{\tau_n\}$ is an admissible sequence of locations. Then there exists a constant $\mu > 0$ such that*

$$\mathcal{K}(F_{\tau_n}, P_{\widehat{\theta}_n}) = O_{\mathbf{P}_F}\left(\frac{\log n}{n(1 - F(\tau_n))}\right) \qquad \text{as } n \to \infty.$$

PROOF. The proof is similar to that of Theorem 4.4. The only changes are that $\widehat{\theta}_n$ replaces $\widehat{\theta}_{n,\tau_n}$ and that one uses (4.15) instead of (4.5). □

## 5. Computation of the rates of convergence.
In this section we shall compute explicitly optimal rates of convergence in two particular cases.

Introduce the distance $\rho_*(x, y) = \max\{|\log \frac{x}{y}|, |\frac{1}{x} - \frac{1}{y}|\}$, $x, y > 0$. From Proposition 8.6 it follows that the sequence $\{\tau_n\}$ is admissible if there exists a function $t \to \theta_t$ such that

$$(5.1) \quad \rho_{\tau_n}^2 \equiv \sup_{x \geq \tau_n} \rho_*(\alpha_F(x), \theta_{\tau_n})^2 = O\left(\frac{\log n}{n(1 - F(\tau_n))}\right) \to 0 \qquad \text{as } n \to \infty,$$

$$(5.2) \qquad \sup_{m \geq n} \int_1^\infty (1 + \log x)^2 x^{r_0} F_{\tau_m}(dx) = O(1) \qquad \text{as } n \to \infty.$$

In turn this implies that the conclusions of Section 4 hold true. The optimal rate corresponds to minimal location $\tau_n$ satisfying (5.1).

### 5.1. *Perturbed Pareto model.*
Assume that $F$ has the form

$$(5.3) \qquad F(x) = 1 - c_\beta x^{-1/\beta} \log x, \qquad x \geq x_0 \geq e,$$

where $\beta \geq \beta_0 > 0$, $x_0$ and $c_\beta$ are chosen such that $F(x)$ is strictly monotone and $F(x_0) = 0$. By straightforward calculations $\theta_t(F) = \beta(1 + \frac{\beta}{\log t})$ and $\alpha_F(x) = \beta(1 - \frac{\beta}{\log x})^{-1}$. Since $\rho_{\tau_n} \leq \frac{\beta}{\log \tau_n}$ and $1 - F(t) = c_\beta t^{-1/\beta} \log t$, for determining nearly optimal location $\tau_n$ we get the balance condition

$$\frac{\beta^2}{\log^2 \tau_n} = O\left(\frac{\log n}{n c_\beta \tau_n^{-1/\beta} \log \tau_n}\right).$$

With the optimal choice $\tau_n \asymp \frac{n^\beta}{\log^3 n}$, one gets $\frac{\log n}{n(1-F(\tau_n))} = O(\log^{-2} n)$ as $n \to \infty$. On the other hand condition (5.2) is satisfied since, by (2.3), $f_{F_{\tau_n}}(x) \leq \frac{1}{\beta} x^{-(1-\varepsilon)/\beta-1}$, for some $\varepsilon \in (0,1)$. Then according to the results in Section 4, $\mathcal{K}(\theta_{\tau_n}(F), \widehat{\theta}_{n,\tau_n})$ and $\mathcal{K}(\theta_{\tau_n}(F), \widehat{\theta}_n)$ are $O_{\mathbf{P}_F}(\log^{-2} n)$ as $n \to \infty$. Taking $\theta_t = \beta$, in the same way one shows that $\mathcal{K}(\beta, \widehat{\theta}_{n,\tau_n})$ and $\mathcal{K}(\beta, \widehat{\theta}_n)$ are $O_{\mathbf{P}_F}(\log^{-2} n)$ as $n \to \infty$. According to the results in Drees [6], Theorem 2.1, the best achievable rate for estimating $\beta$ in $L^2$ norm is $\frac{1}{\log n}$, in a certain class of d.f.'s which includes the d.f. $F$ satisfying (5.3) (we refer to Drees [6] for details). Our estimators $\widehat{\theta}_{n,\tau_n}$ and $\widehat{\theta}_n$ attain the same rate.

For the log-gamma d.f. we obtain the same rate of convergence since it has essentially the same behavior as the d.f. $F$ defined by (5.3).

5.2. *Hall model.* Assume that $F$ is of the form

$$(5.4) \qquad F(x) = 1 - c_\beta x^{-1/\beta} - c_\gamma x^{-1/\gamma}, \qquad x \geq x_0,$$

where $\beta = \gamma + \alpha \geq \beta_0 > 0$, $\alpha, \beta, \gamma > 0$ and $x_0$, $c_\beta$ and $c_\gamma$ are such that $F(x)$ is increasing on $[x_0, \infty)$. Also, though it is not exactly the model proposed by Hall, we shall call it the Hall model. By straightforward calculations $\theta_t(F) = \frac{\beta c_\beta t^{-1/\beta} + \gamma c_\gamma t^{-1/\gamma}}{c_\beta t^{-1/\beta} + c_\gamma t^{-1/\gamma}}$ and $\alpha_F(x) = \frac{c_\beta x^{-1/\beta} + c_\gamma x^{-1/\gamma}}{\beta^{-1} c_\beta x^{-1/\beta} + \gamma^{-1} c_\gamma x^{-1/\gamma}}$. It is easy to check that $\rho_{\tau_n} = O(\tau_n^{-1/\gamma+1/\beta})$ as $n \to \infty$. Since $1 - F(t) = c_\beta t^{-1/\beta} + c_\gamma t^{-1/\gamma}$, for determining the nearly optimal location $\tau_n$ we get the balance condition

$$\tau_n^{-2/\gamma+2/\beta} = O\left(\frac{\log n}{n(c_\beta \tau_n^{-1/\beta} + c_\gamma \tau_n^{-1/\gamma})}\right).$$

The optimal choice $\tau_n \asymp (\frac{n}{\log n})^{\beta\gamma/(2\beta-\gamma)}$, implies $\frac{\log n}{n(1-F(\tau_n))} = O((\frac{\log n}{n})^{2\alpha/(\beta+\alpha)})$ as $n \to \infty$. As in the previous example one can show that (5.2) is satisfied. Then according to the results in Section 4, $\mathcal{K}(\theta_{\tau_n}(F), \widehat{\theta}_{n,\tau_n})$ and $\mathcal{K}(\theta_{\tau_n}(F), \widehat{\theta}_n)$ are $O_{\mathbf{P}_F}((\frac{\log n}{n})^{2\alpha/(\beta+\alpha)})$ as $n \to \infty$. Taking $\theta_t = \beta$ we have that $\mathcal{K}(\beta, \widehat{\theta}_{n,\tau_n})$ and $\mathcal{K}(\beta, \widehat{\theta}_n)$ are $O_{\mathbf{P}_F}((\frac{\log n}{n})^{2\alpha/(\beta+\alpha)})$ as $n \to \infty$. By the results in Hall and Welsh [12] (see also Drees [6], Theorem 2.1, for a more general result), the optimal rate of convergence that can be achieved for estimating $\alpha = 1/\beta$ in $L^2$ norm is $n^{\rho/(2\rho+1)}$, in the class of d.f. $F$ having the density (3.1). The d.f. $F$ defined by (5.4) satisfies this condition with $\gamma = \beta(1+\rho)^{-1}$. Since $\frac{\alpha}{\beta+\alpha} = \frac{2\beta-2\gamma}{2\beta-\gamma} = \frac{\rho}{2\rho+1}$, the estimators $\widehat{\theta}_{n,\tau_n}$ and $\widehat{\theta}_n$ attain this rate for $\beta$ up to an additional $\log^{\alpha/(\beta+\alpha)} n$ factor.

## 6. Numerical results.

6.1. *Choice of the parameters of the adaptive procedure.* An important parameter in the proposed adaptive procedure is the sequence of critical values $\mathfrak{z}_n$. According to Remark 3.1 the test statistic does not depend on the parameter of the Pareto law if the observations follow a Pareto model. Therefore we propose to compute the critical values $\mathfrak{z}_n$ by Monte Carlo simulations from the homogeneous model with i.i.d. standard Pareto observations $X_i$, $i = 1, \ldots, n$.

We simulated 2000 realizations with three different sample sizes $n = 200$, 500, 1000 and with the grid length $K_n$ set to 200. The size of testing windows in $T_{n,m}$ is determined by $\rho = 1/4$ and $\delta = 1/20$. The empirical d.f. of the statistic $T_n = \max_{m=r_1, \ldots, r_{K_k}} T_{n,m}$ has been computed and it was found that in all simulations the critical value $\mathfrak{z}_n = 10$ corresponds to a 99% confidence level. The same critical value $\mathfrak{z}_n = 10$, corresponding to a 99% confidence level, has been found from 2000 realizations with $n = 1000$, $\rho = 1/4$, $\delta = 1/20$ and with different grid lengths $K_n = 100$, 200, 300. Additional simulations show that finite sample properties of the test statistic $T_n$ do depend very little on the parameters $k_0$ and $\delta$. The value $\mathfrak{z}_n = 10$ which approximately corresponds to a 99% confidence level in all cases and the grid length $K_n = 200$ have been retained.

Further simulations show that the finite sample performance of the adaptive estimator depends mainly on the parameter $\rho$ which plays the same role as the bandwidth in the nonparametric kernel density estimation. The choice of $\rho$, in turn, depends on the class of functions in hands. In the simulations below we fix the following values $\delta = 1/20$, $k_0 = n/20$, $K_n = 200$, $\mathfrak{z}_n = 10$. As to the value of $\rho$ it will be fixed to $1/4$. This choice is motivated by the desire to minimize the relative mean squared error for some given heavy-tailed laws. In the simulations below we shall consider the following distributions: (1) The positive part of Cauchy d.f. $F(x) = \frac{2}{\pi} \arctan x$, $x \geq 0$. (2) Log-gamma d.f. $F(x) = G_{1,2}(\log x)$, $x \geq 1$, where $G_{\lambda,\alpha}(x)$, $x \geq 0$, is gamma d.f. with parameters $\lambda, \alpha > 0$. (3) Log perturbed Pareto d.f. $F(x) = 1 - x^{-1} \log x$, $x \geq x_0 = e$. (4) Hall's model $F(x) = 1 - 2x^{-1} + x^{-2.5}$, $x \geq x_0 > 0$, where $x_0$ satisfies $2x_0^{-1} - x_0^{-2.5} = 1$. (5) GPD $F(x) = 1 - (1+x)^{-1}$, $x \geq 0$.

6.2. *Estimation of extreme quantiles.* We shall demonstrate the performance of the adaptive estimator $\widehat{\theta}_n = \widehat{h}_{n,\widehat{k}_n}$ by presenting the results of a simulation study for estimating extreme quantiles. We consider two opposite cases: observations from d.f.'s whose tails are close to a Pareto model in the range of big order statistics (such as Cauchy d.f., GPD, some of the Hall models) and observations from d.f.'s whose tails are not well approximated by a Pareto model in the range of the large order statistics at least for samples of reasonable size (such as log-gamma d.f. and log perturbed Pareto d.f.). For many d.f.'s our simulations show a behavior in-between the latter two types. We performed 2000 Monte Carlo simulations of $n = 1000$ observations. The quantiles of $F$ are estimated by solving for $x \geq \tau_n$ in the following approximation formula:

$$1 - P_{\widehat{\theta}_n}\left(\frac{x}{\tau_n}\right) = \left(\frac{x}{\tau_n}\right)^{-1/\widehat{\theta}_n} \approx 1 - F_{\tau_n}\left(\frac{x}{\tau_n}\right) = \frac{1 - F(x)}{1 - F(\tau_n)} = \frac{1 - p}{1 - F(\tau_n)}.$$

If $x < \tau_n$ we determine $x$ from the equality $p = F(x)$. The unknown location parameter $\tau_n$ has to be replaced with the adaptive value $\widehat{\tau}_n = X_{n,\widehat{k}_n}$ and $F$ with the empirical d.f. $\widehat{F}_n$, which leads to the following adaptive estimate of the quantiles of $F$:

$$(6.1) \quad \widehat{q}_{n,p} = \widehat{q}_{n,\widehat{k}_n,p} = \begin{cases} X_{n,[n(1-p)]}, & \text{if } p < 1 - \dfrac{\widehat{k}_n}{n}, \\ X_{n,\widehat{k}_n}\left(\dfrac{\widehat{k}_n}{n(1-p)}\right)^{\widehat{\theta}_n}, & \text{otherwise,} \end{cases} \quad p \in (0,1).$$

Here and in the sequel $\widehat{q}_{n,k,p}$ denotes the quantile estimator

$$\widehat{q}_{n,k,p} = \begin{cases} X_{n,[n(1-p)]}, & \text{if } p < 1 - \dfrac{k}{n}, \\ X_{n,k}\left(\dfrac{k}{n(1-p)}\right)^{\widehat{h}_{n,k}}, & \text{otherwise,} \end{cases} \quad p \in (0,1), k = 2, \ldots, n,$$

(6.2)

which combines the sample quantile estimator for low quantiles and the estimator introduced by Weissman [21] for high quantiles.

6.2.1. *The performance of the adaptive estimator.* For any estimator $\widehat{\alpha}$ of $\alpha$ let

$$\sigma^2(\widehat{\alpha}, \alpha) = \frac{1}{n} \sum_{i=1}^{n} \log^2 \frac{\widehat{\alpha}}{\alpha}$$

be the relative mean squared error (RelMSE) of $\widehat{\alpha}$. We compare $\sigma(\widehat{q}_{n,p}, q_p)$ with $\sigma(\widehat{q}_{n,k,p}, q_p)$. Figures 4 and 5 plot these quantities for $p = 1 - 1/n = 0.999$ and $p = 0.9999999$ as a function of $k$. It is useful to compare RelMSE $\sigma(\widehat{q}_{n,p}, q_p)$ with minimal RelMSE $\min_k \sigma(\widehat{q}_{n,k,p}, q_p)$ as a function of $p$ (see Figure 6). The ratio $r_{n,p} = \sigma(\widehat{q}_{n,p}, q_p) / \min_k \sigma(\widehat{q}_{n,k,p}, q_p)$ regarded as a function of $p$ is plotted in Figure 7 (see also Table 1 for a more precise evaluation). These simulations show that the proposed adaptive procedure captures nearly the best choice in $k$ which depends on the unknown d.f. $F$. The procedure gives reasonable results in both cases, for d.f. with Pareto like tails as well as with d.f. which exhibits large perturbations from these tails. Table 1 hints that the increase of RelMSE introduced by the adaptive procedure for estimating high quantiles $q_p$ with $p \in [0.9, 0.9999999999]$ does not exceed 7%. For log-perturbed Pareto d.f. the results are very similar to those of log-gamma d.f. and therefore will not be presented here.

We would like to point out that for GPD the ratio $r_{n,p}$ is even less than 1, which means that the adaptive quantile estimator $\hat{q}_{n,p}$ improves the performance of individual quantile estimators $\hat{q}_{n,k,p}$, $k = 2, \ldots, n$. This improvement can be observed for other d.f. with an appropriate choice of the parameter $\rho$. The corresponding plots of the ratio $r_{n,p}$ for log-gamma d.f. with $\rho = 1/10$ and for Hall model with $\rho = 1/2$ are given in Figure 8.

The high variability for extreme quantiles $q_p$, $p > 1 - 0.1/n$ (see Figure 6) is mainly explained by the bias introduced by the Pareto model and less by the variability introduced by the adaptive procedure. The bias reducing techniques can be applied under some additional assumptions on the underlying d.f. $F$. Our adaptive values $\hat{k}_n$ and $\hat{\theta}_n$ can be applied with these types of bias reduced estimators to construct new adaptive quantile estimators, however this issue will not be discussed here. For further details on this subject we refer to Danielsson et al. [3], Gomes and Oliveira [8]; see also Chapter 4.7 in Beirlant et al. [1].

6.2.2. *Comparison with sample quantiles.* For any $k = 1, \ldots, n$ the sample quantile $X_{n,k}$ is considered as an estimate of the true quantile $F^{-1}(p_{n,k})$, where $p_{n,k} = 1 - k/n$. We shall compare the RelMSE of adaptive quantiles $\hat{q}_{n,p_{n,k}}$, with those of sample quantiles $X_{n,k}$, for $k = 1, \ldots, 500$ by computing the ratio $r^0_{n,k} = \sigma(X_{n,k}, q_{p_{n,k}})/\sigma(\hat{q}_{n,p_{n,k}}, q_{p_{n,k}})$. The results of the simulations are reported in Table 2 and Figure 9. They show that there is a substantial gain in variance if we use (6.1) for estimating large quantiles.

Figures reported in Tables 1 and 2 can be used to compare the performance of the adaptive estimator $\hat{\theta}_n$ with other adaptive estimators.

TABLE 1
*Values of $r_{n,p}$*

| $p$ | 0.9 | 0.99 | 0.999 | 0.9999 | 0.99999 |
|---|---|---|---|---|---|
| Cauchy | 1.017966 | 1.023952 | 1.041944 | 1.049905 | 1.054291 |
| log-gamma | 1.042706 | 1.002527 | 1.002542 | 1.013393 | 1.021253 |
| Hall model | 0.996002 | 1.009698 | 1.023196 | 1.030144 | 1.034276 |
| GPD | 1.094321 | 0.998349 | 0.989391 | 0.985767 | 0.984071 |
| $p$ | 0.999999 | 0.9999999 | 0.99999999 | 0.999999999 | 0.9999999999 |
| Cauchy | 1.057159 | 1.059174 | 1.060642 | 1.061758 | 1.062635 |
| log-gamma | 1.026952 | 1.031355 | 1.03472 | 1.037275 | 1.039637 |
| Hall model | 1.036994 | 1.038913 | 1.040339 | 1.041438 | 1.042312 |
| GPD | 0.983118 | 0.982513 | 0.982184 | 0.981981 | 0.981829 |

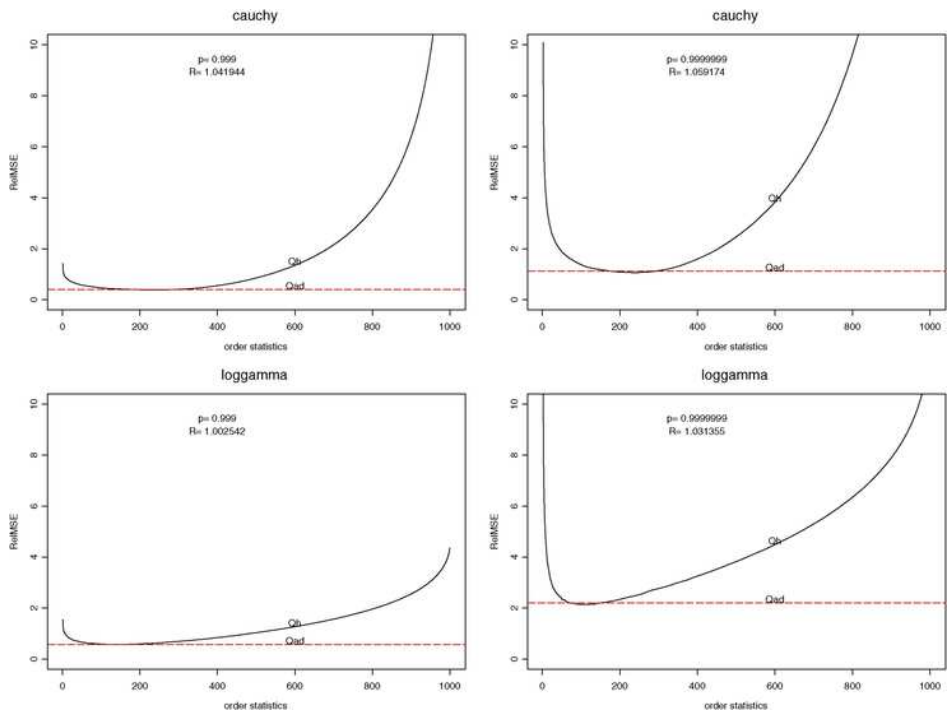

Fig. 4. $Qh = \sigma(\widehat{q}_{n,k,p}, q_p)$ for $k = 1, \ldots, n$ and $Qad = \sigma(\widehat{q}_{n,p}, q_p)$; $R = r_{n,p}$. Top: Cauchy observations; Bottom: log-gamma observations.

6.3. *Estimation of the index of regular variation.* According to (4.16) the adaptive estimator $\widehat{\theta}_n$ converges to the index of regular variation $\gamma$. The performance of an estimator $\widetilde{\theta}_n$ w.r.t. $\gamma$ will be measured using the root

Table 2
*Values of $r_{n,k}^0$*

| $k$ | 1 | 2 | 3 | 4 | 5 | 10 | 20 |
|---|---|---|---|---|---|---|---|
| Cauchy | 3.5360 | 2.4100 | 2.0294 | 1.8671 | 1.7200 | 1.4226 | 1.2621 |
| log-gamma | 2.7270 | 1.9417 | 1.7306 | 1.5924 | 1.4971 | 1.3010 | 1.2453 |
| Hall model | 4.1240 | 2.7809 | 2.3237 | 2.1246 | 1.9466 | 1.5772 | 1.3605 |
| GPD | 1.3117 | 1.2108 | 2.9563 | 2.0609 | 1.7629 | 1.6422 | 1.5288 |
| $k$ | 30 | 40 | 50 | 60 | 70 | 80 | 90 |
| Cauchy | 1.2081 | 1.1852 | 1.1849 | 1.1755 | 1.1928 | 1.1860 | 1.1745 |
| log-gamma | 1.1982 | 1.1724 | 1.1611 | 1.1696 | 1.1642 | 1.1622 | 1.1748 |
| Hall model | 1.2758 | 1.2324 | 1.2183 | 1.2001 | 1.2141 | 1.2088 | 1.2040 |
| GPD | 1.1813 | 1.1683 | 1.1675 | 1.1509 | 1.1557 | 1.1327 | 1.1033 |

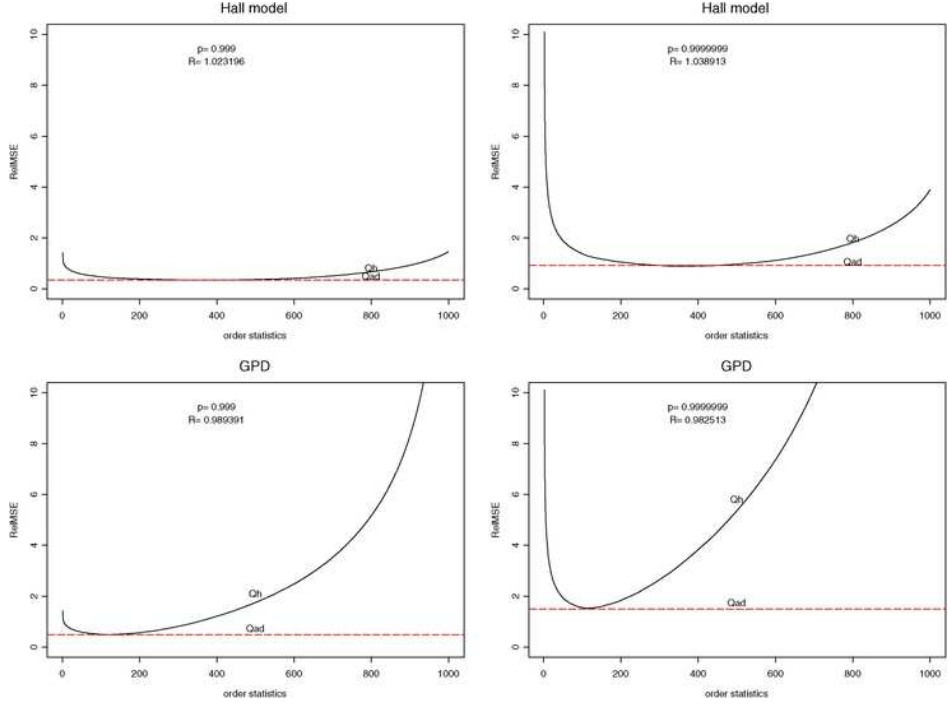

Fig. 5.   $Qh = \sigma(\widehat{q}_{n,k,p}, q_p)$ for $k = 1, \ldots, n$ and $Qad = \sigma(\widehat{q}_{n,p}, q_p)$; $R = r_{n,p}$. Top: Observations from Hall's model; Bottom: GPD observations.

mean squared error (RMSE) $\sigma(\widetilde{\theta}_n) = E^{1/2}(\widetilde{\theta}_n - \gamma)^2$. The corresponding simulations of the RMSE's $\sigma(\widehat{\theta}_n)$ and $\sigma(\widehat{h}_{n,k})$ (as a function of $k$) are presented in Figure 10. In case of Cauchy d.f. the minimal value of RMSE of the Hill estimator is $\min_k \sigma(\widehat{h}_{n,k}) = 0.07385$, while the RMSE of the adaptive estimator is $\sigma(\widehat{\theta}_n) = 0.07899$, which gives the ratio $r_n^\gamma = \sigma(\widehat{\theta}_n)/\min_k \sigma(\widehat{h}_{n,k}) =$

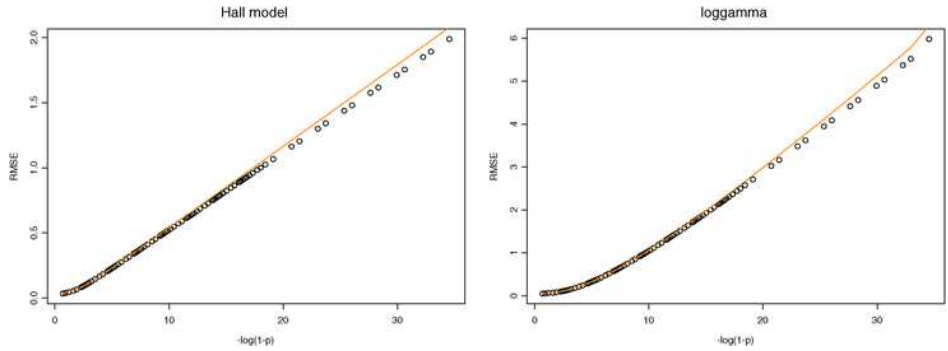

Fig. 6.   $\min_k \sigma(\widehat{q}_{n,k,p}, q_p)$ (points) and $\sigma(\widehat{q}_{n,p}, q_p)$ (solid line) as functions of $p$.

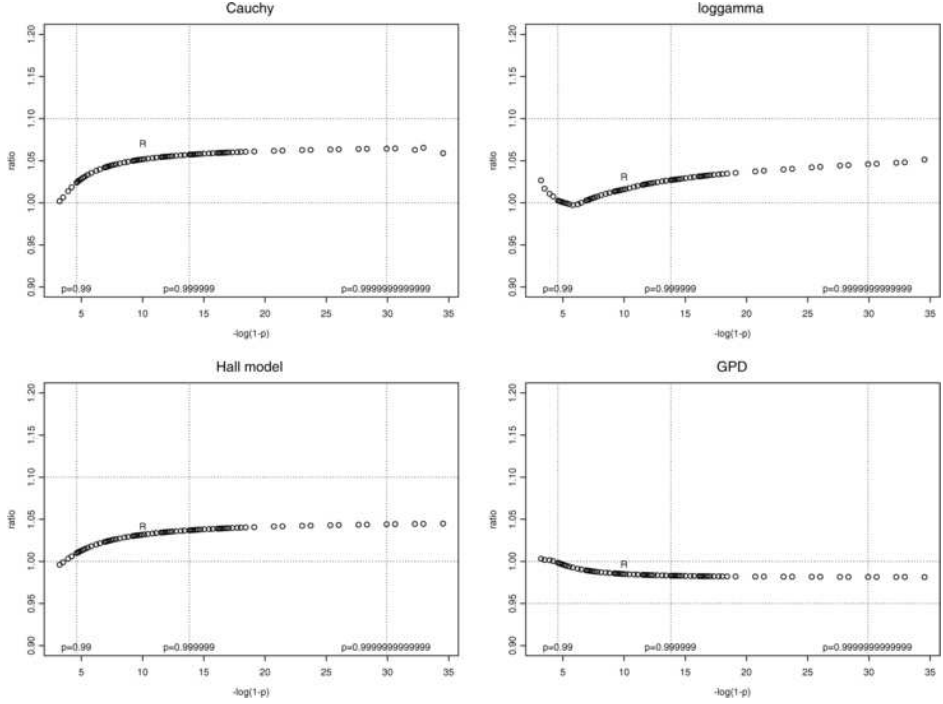

FIG. 7.    *The ratio $r_{n,p}$ as a function of p.*

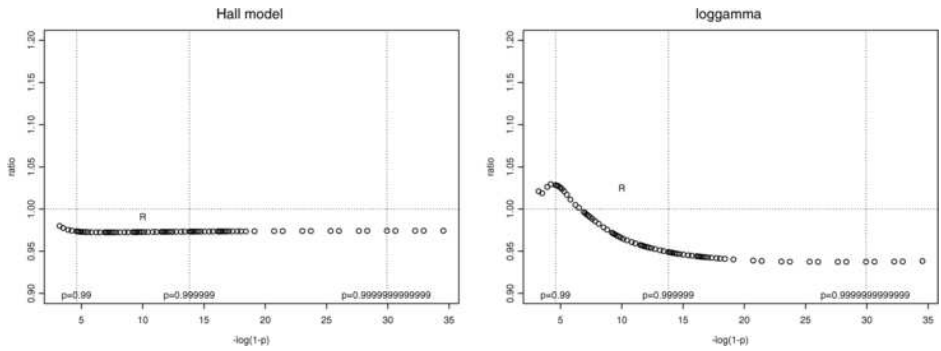

FIG. 8.    *The ratio $r_{n,p}$ as a function of p. Adaptive procedure is performed with $\rho = 1/2$ (left) and $\rho = 1/10$ (right).*

1.06966. For log-gamma d.f. the minimal value of RMSE of the Hill estimator is $\min_k \sigma(\widehat{h}_{n,k}) = 0.23112$, while the RMSE of the adaptive estimator is $\sigma(\widehat{\theta}_n) = 0.24804$, which gives the ratio $r_n^\gamma = \sigma(\widehat{\theta}_n)/\min_k \sigma(\widehat{h}_{n,k}) = 1.07321$. Thus for Cauchy and log-gamma the adaptive estimator increases the minimal variance in the family of Hill estimators by 7.4%.

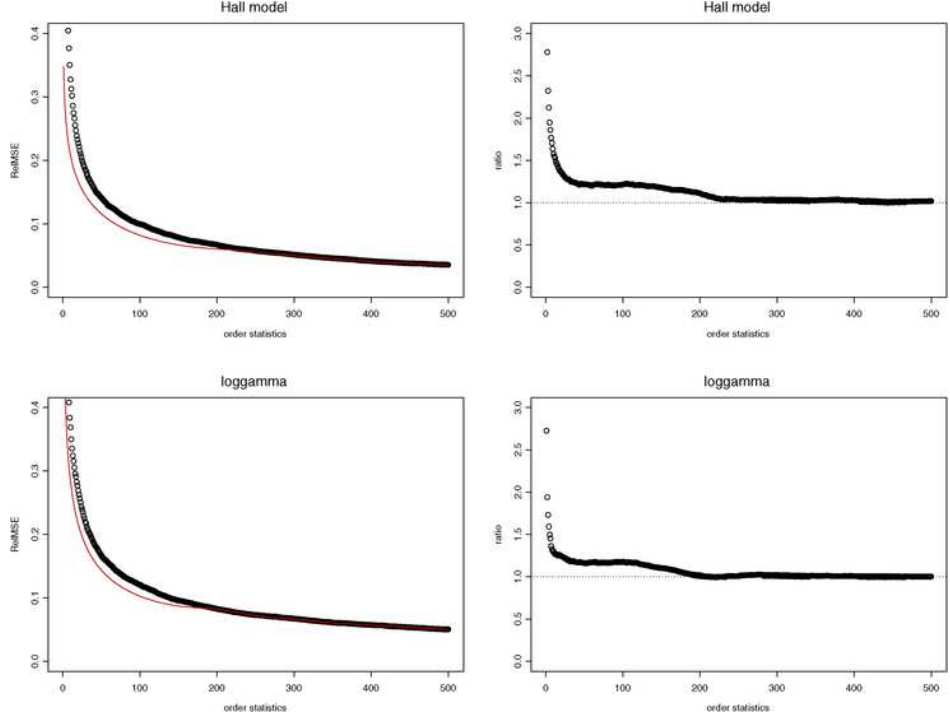

Fig. 9.   *Left: $\sigma(X_{n,k}, q_{p_{n,k}})$ (points) and $\sigma(\widehat{q}_{n,p_{n,k}}, q_{p_{n,k}})$ (solid line) as functions of k; Right: $r_{n,k}^0 = \sigma(X_{n,k}, q_{p_{n,k}})/\sigma(\widehat{q}_{n,p_{n,k}}, q_{p_{n,k}})$ as function of k; Top: Hall model; Bottom: log-gamma d.f.*

**7. Proofs of the exponential bounds.**   Let $t \geq x_0$. The local log-likelihood ratio $L_{n,t}(H, G) = L_{n,t}(H) - L_{n,t}(G)$ admits the representation

$$L_{n,t}(H, G) = \sum_{i:X_i > t} \log \frac{\alpha_G(X_i)}{\alpha_H(X_i)} + \int_{(t, X_i]} \left( \frac{1}{\alpha_G(u)} - \frac{1}{\alpha_H(u)} \right) \frac{du}{u}.$$

Recall the following notations: $n_t = n(1 - F(t))$, $\widehat{n}_t = \sum_{i=1}^n 1(X_i > t)$ and $\widehat{n}_{t,\tau} = \sum_{i=1}^n 1(t < X_i \leq \tau)$ for $t \geq x_0$.

We start with a bound for the log of the local likelihood ratio.

PROPOSITION 7.1.   *Let $s \geq x_0$. For any $F, G, H \in \mathcal{F}$ any $y > 0$ it holds,*

$$(7.1) \qquad \mathbf{P}_F(L_{n,s}(H, G) > y) \leq \exp\left( -\frac{y}{2} + \frac{n_s}{2} d_s \right),$$

*where $d_s = \chi^2(F_s, G_s)$.*

PROOF.   By exponential Chebyshev's inequality,

$$\mathbf{P}_F(L_{n,s}(H, G) > y) \leq \exp(-y/2 + \log \mathbf{E}_F(\exp(L_{n,s}(H, G)/2))).$$

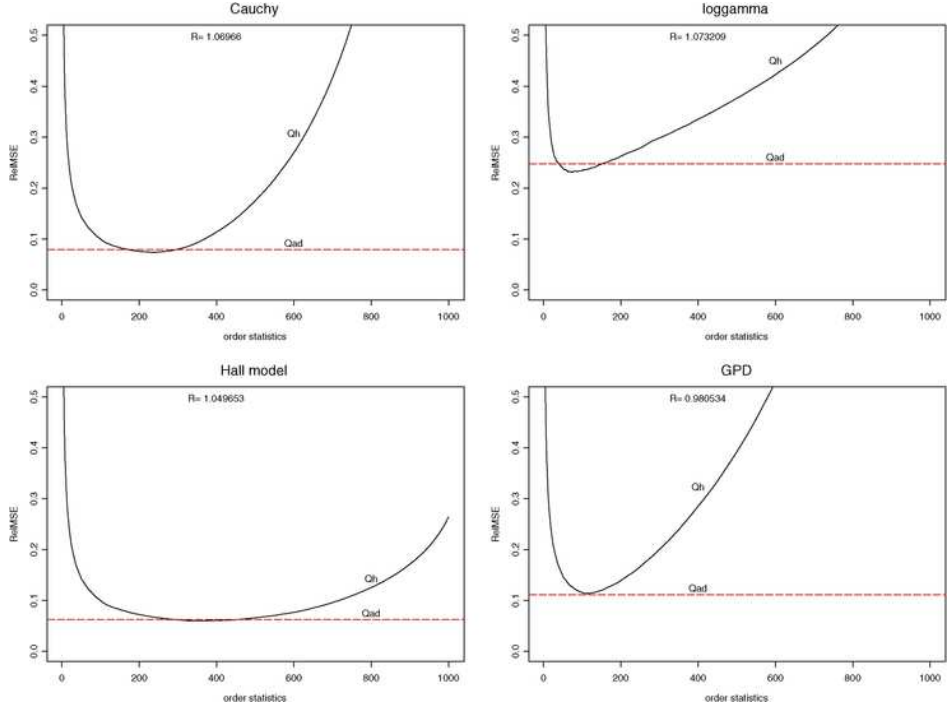

Fig. 10. *Estimation of the index of regular variation $\gamma$: $Qh = \sigma(\widehat{h}_{n,k})$—plot of the estimated RMSE of the Hill estimator $\widehat{h}_{n,k}$, for $k = 1, \ldots, n$; $Qad = \sigma(\widehat{\theta}_n)$—RMSE of the adaptive estimator $\widehat{\theta}_n$ w.r.t. $\gamma$.*

Since the r.v.'s $X_1, \ldots, X_n$ are i.i.d., one gets

$$\log \mathbf{E}_F\left(\exp\left(\frac{1}{2} L_{n,s}(H, G)\right)\right) = n \log \mathbf{E}_F\left(\exp\left(\frac{1_{A_{i,s}}}{2} \log \frac{dH_s}{dG_s}\left(\frac{X_i}{s}\right)\right)\right),$$

where $A_{i,s} = \{X_i > s\}$. Since $\mathbf{E}_F \exp(1_{A_{i,s}} \log \frac{dH_s}{dF_s}(\frac{X_i}{s})) = 1$, by Hölder's inequality

$$\mathbf{E}_F \exp\left(\frac{1_{A_{i,s}}}{2} \log \frac{dH_s}{dG_s}\left(\frac{X_i}{s}\right)\right) \leq \mathbf{E}_F^{1/2} \exp\left(1_{A_{i,s}} \log \frac{dF_s}{dG_s}\left(\frac{X_i}{s}\right)\right).$$

Using

$$\mathbf{E}_F \exp\left(1_{A_{i,s}} \log \frac{dF_s}{dG_s}\left(\frac{X_i}{s}\right)\right) = F(s) + (1 - F(s))(d_s + 1),$$

one gets

$$\mathbf{P}_F(L_{n,s}(H, G) > y) \leq \exp\left(-\frac{y}{2} + \frac{n}{2} \log\{1 + (1 - F(s))d_s\}\right)$$

$$\leq \exp\left(-\frac{y}{2} + \frac{n_s}{2} d_s\right).$$

$\square$

Denote for brevity $L_{n,t,\tau}(\theta', \theta) = L_{n,t}(\theta', \theta, \tau, \theta)$.

COROLLARY 7.2.    *For any $\tau \geq t \geq s \geq x_0$ and $\theta', \theta > 0$,*

$$(7.2) \qquad \mathbf{P}_F(L_{n,t}(\theta', \theta) > y) \leq \exp\left(-\frac{y}{2} + \frac{n_s}{2} d_s\right),$$

$$(7.3) \qquad \mathbf{P}_F(L_{n,t,\tau}(\theta', \theta) > y) \leq \exp\left(-\frac{y}{2} + \frac{n_s}{2} d_s\right),$$

*where $d_s = \chi^2(F_s, P_\theta)$.*

PROOF.    The first assertion follows from (7.1) when applied with $H$ and $G$ such that $\alpha_H(x) = \alpha_G(x) = \theta$ for $x \in [x_0, t)$ and $\alpha_H(x) = \theta'$, $\alpha_G(x) = \theta$ for $x \in [t, \infty)$. The second one is obtained with $\alpha_H(x) = \alpha_G(x) = \theta$ for $x \in [x_0, t) \cup [\tau, \infty)$ and $\alpha_H(x) = \theta'$, $\alpha_G(x) = \theta$ for $x \in [t, \tau)$.  $\square$

PROPOSITION 7.3.    *For any $F \in \mathcal{F}$, $\theta > 0$, $y > 0$, $\tau \geq t \geq s \geq x_0$ it holds*

$$\mathbf{P}_F(\widehat{n}_t \mathcal{K}(\widehat{\theta}_{n,t}, \theta) > y) \leq 2n \exp\left(-\frac{y}{2} + \frac{n_s}{2} d_s\right),$$

$$\mathbf{P}_F(\widehat{n}_{t,\tau} \mathcal{K}(\widehat{\theta}_{n,t,\tau}, \theta) > y) \leq 2n \exp\left(-\frac{y}{2} + \frac{n_s}{2} d_s\right),$$

*where $d_s = \chi^2(F_s, P_\theta)$.*

PROOF.    We shall prove only the second inequality the first one being proved in the same way.

First note that $\widehat{n}_{t,\tau} \mathcal{K}(\widehat{\theta}_{n,t}, \theta) = L_{n,t,\tau}(\widehat{\theta}_{n,t,\tau}, \theta)$. For the sake of brevity let $l_k(\alpha) = (y/k - \log\frac{\theta}{\alpha})/(\frac{1}{\theta} - \frac{1}{\alpha})$, $\alpha > \theta$. Since the function $l_k(\alpha)$ is continuous in $\alpha$ for $\alpha > \theta$ and $\lim_{\alpha \to \theta} l_k(\alpha) = \infty$, $\lim_{\alpha \to \infty} l_k(\alpha) = \infty$, there exists a finite point $\alpha_k^* > \theta$ which realize $\alpha_k^* = \arg\min_{\alpha > \theta} l_k(\alpha)$. Note that $\alpha_k^*$ is a function only on $k, y, \theta$. With these notation, on the event $A_k = \{\widehat{\theta}_{n,t,\tau} > \theta, \ \widehat{n}_{t,\tau} = k\}$, we have that the inequality

$$L_{n,t,\tau}(\widehat{\theta}_{n,t,\tau}, \theta) = \widehat{n}_{t,\tau}\left(\log\frac{\theta}{\widehat{\theta}_{n,t,\tau}} - (1/\widehat{\theta}_{n,t,\tau} - 1/\theta)\widehat{\theta}_{n,t,\tau}\right) > y$$

is equivalent to $\widehat{\theta}_{n,t,\tau} \geq l_k(\widehat{\theta}_{n,t,\tau})$ and the inequality $\widehat{\theta}_{n,t,\tau} \geq l_k(\alpha_k^*)$ is equivalent to $L_{n,t,\tau}(\alpha_k^*, \theta) > y$. Then

$$\{L_{n,t,\tau}(\widehat{\theta}_{n,t,\tau}, \theta) > y\} \cap A_k = \{\widehat{\theta}_{n,t,\tau} \geq l_k(\widehat{\theta}_{n,t,\tau})\} \cap A_k$$

$$\subseteq \{\widehat{\theta}_{n,t,\tau} \geq l_k(\alpha_k^*)\} \cap A_k$$

$$\subseteq \{L_{n,t,\tau}(\alpha_k^*, \theta) > y\}.$$

In the same way

$$\{L_{n,t,\tau}(\widehat{\theta}_{n,t,\tau}, \theta) > y\} \cap B_k \subseteq \{L_{n,t,\tau}(\alpha_k^{**}, \theta) > y\},$$

where $B_k = \{\widehat{\theta}_{n,t,\tau} \leq \theta, \ \widehat{n}_{t,\tau} = k\}$ and $\alpha_k^{**} = \arg\max_{0 < \alpha \leq \theta} l_k(\alpha)$ is a function only on $k, y, \theta$. The latter implies

$$\mathbf{P}_F(L_{n,t,\tau}(\widehat{\theta}_{n,t,\tau}, \theta) > y, \ \widehat{n}_{t,\tau} = k)$$
$$\leq \mathbf{P}_F(L_{n,t,\tau}(\alpha_k^*, \theta) > y) + \mathbf{P}_F(L_{n,t,\tau}(\alpha_k^{**}, \theta) > y).$$

Since by Corollary 7.2,

$$\mathbf{P}_F(L_{n,t,\tau}(\theta', \theta) > x) \leq \exp\left(-\frac{y}{2} + \frac{n_s}{2} d_s\right),$$

with $\theta' = \alpha_k^*, \alpha_k^{**}$, one gets

$$\mathbf{P}_F(L_{n,t,\tau}(\widehat{\theta}_{n,t,\tau}, \theta) > x) = \sum_{k=1}^{n} \mathbf{P}_F(L_{n,t,\tau}(\widehat{\theta}_{n,t,\tau}, \theta) > x, \ \widehat{n}_{t,\tau} = k)$$
$$\leq 2n \exp\left(-\frac{y}{2} + \frac{n_s}{2} d_s\right),$$

which completes the proof.  $\square$

PROPOSITION 7.4. *For any $F \in \mathcal{F}$ and $s \geq x_0$, $\theta > 0$, $y > 0$ it holds*

$$\mathbf{P}_F\left(\sup_{s \leq t} \widehat{n}_t \mathcal{K}(\widehat{\theta}_{n,t}, \theta) > y\right) \leq 2n^4 \exp\left(-\frac{y}{2} + \frac{n_s}{2} d_s\right) + \frac{1}{n},$$

$$\mathbf{P}_F\left(\sup_{s \leq t \leq \tau} \widehat{n}_{t,\tau} \mathcal{K}(\widehat{\theta}_{n,t,\tau}, \theta) > y\right) \leq n^7 \exp\left(-\frac{y}{2} + \frac{n_s}{2} d_s\right) + \frac{1}{n},$$

*where $d_s = \chi^2(F_s, P_\theta)$.*

PROOF. We shall give a proof only for the second inequality, the first one being proved in the same way.

Let $N = n^3$ and $J = \{s_0, \dots, s_N\}$ be the set of numbers satisfying $s_{i-1} < s_i$, $F([s_{i-1}, s_i)) = 1/N$ and $\bigcup_{i=1}^{N} [s_{i-1}, s_i) = [x_0, \infty)$. If we denote by $\mathfrak{A}_n$ the event that $X_{n,1}, \dots, X_{n,n}$ will fall into disjoint intervals, then, for $K_n > 2$,

$$\mathbf{P}_F(\mathfrak{A}_n) = \prod_{i=1}^{n}\left(1 - \frac{i-1}{N}\right) \geq 1 - \sum_{i=1}^{n} \log\left(1 - \frac{i-1}{N}\right)$$
$$\geq 1 - \frac{3}{2}\sum_{i=2}^{n} \frac{i-1}{N} = 1 - \frac{3n(n-1)}{4n^3} \geq 1 - \frac{1}{n}.$$

On the event $\mathfrak{A}_n$ it holds

$$\sup_{s \leq t \leq \tau} \widehat{n}_{t,\tau} \mathcal{K}(\widehat{\theta}_{n,t,\tau}, \theta) = \max_{s \leq t \leq \tau, \ t, \tau \in J} \widehat{n}_{t,\tau} \mathcal{K}(\widehat{\theta}_{n,t,\tau}, \theta).$$

Then

$$(7.4) \quad \mathbf{P}_F\left(\sup_{s \le t \le \tau} \widehat{n}_{t,\tau} \mathcal{K}(\widehat{\theta}_{n,t,\tau}, \theta) > y\right)$$

$$\le \sum_{s \le t \le \tau, \ t,\tau \in J} \mathbf{P}_F(\widehat{n}_{t,\tau} \mathcal{K}(\widehat{\theta}_{n,t,\tau}, \theta) > y) + 1 - \mathbf{P}_F(\mathfrak{A}_n).$$

According to Proposition 7.3

$$\mathbf{P}_F(L_{n,t,\tau}(\widehat{\theta}_{n,t,\tau}, \theta) \ge y) \le 2n \exp(-y_s),$$

where $y_s = \frac{y}{2} - \frac{n}{2}(1 - F(s))d_s$. Since $\sum_{s \le t \le \tau, \ t,\tau \in J} \le n^6/2$, from (7.4) one gets

$$\mathbf{P}_F\left(\sup_{s \le t \le \tau} \widehat{n}_{t,\tau} \mathcal{K}(\widehat{\theta}_{n,t,\tau}, \theta) > y\right) \le n^7 \exp(-y_s) + \frac{1}{n}. \qquad \square$$

We end this section with an exponential bound for the statistic $T_n(t, \tau)$.

PROPOSITION 7.5. *For any $F \in \mathcal{F}$ and $s \ge x_0$, $\theta > 0$, $y > 0$ it holds*

$$\mathbf{P}_F\left(\sup_{s \le t \le \tau} T_n(t, \tau) > 2y\right) \le 2n^7 \exp\left(-\frac{y}{2} + \frac{n_s}{2}d_s\right) + \frac{2}{n},$$

*where $d_s = \chi^2(F_s, P_\theta)$.*

PROOF. Let $\theta > 0$. Using (3.2) and the inequality $\sup_{F \in \mathcal{F}_t} L_{n,t}(F) \ge L_{n,t}(\theta)$, one gets $T_n(t, \tau) \le L_{n,s}(\widehat{\theta}_{n,t,\tau}, \widehat{n}_{n,\tau}, \tau) - L_{n,t}(\theta) = L_{n,t}(\widehat{\theta}_{n,t,\tau}, \widehat{\theta}_{n,\tau}, \tau, \theta)$. The representation (2.8) implies $T_n(t, \tau) \le \widehat{n}_{t,\tau} \mathcal{K}(\widehat{\theta}_{n,t,\tau}, \theta) + \widehat{n}_\tau \mathcal{K}(\widehat{\theta}_{n,\tau}, \theta)$. The assertion of the lemma follows from Proposition 7.4. $\square$

## 8. Auxiliary statements.

LEMMA 8.1. *For any $\theta_1, \theta_2 > 0$ such that $\mathcal{K}(\theta_1, \theta_2) \le \frac{1}{2}$ it holds*

$$(8.1) \quad \frac{1}{3}\log^2 \frac{\theta_1}{\theta_2} \le \mathcal{K}(\theta_1, \theta_2)$$

*and for any $\theta_1, \theta_2 > 0$ such that $\log^2 \frac{\theta_1}{\theta_2} \le \frac{2}{3}$, it holds*

$$(8.2) \quad \mathcal{K}(\theta_1, \theta_2) \le \frac{3}{4}\log^2 \frac{\theta_1}{\theta_2}.$$

*In particular, for any $\theta_1, \theta_2 > 0$ such that $\mathcal{K}(\theta_1, \theta_2) \le \frac{1}{2}$ it holds*

$$\mathcal{K}(\theta_1, \theta_2) \le \frac{9}{4}\mathcal{K}(\theta_2, \theta_1).$$

PROOF. Note that $\frac{1}{3}\log^2(x+1) \le G(x)$ for any $x$ satisfying $G(x) \le 1/2$ and $G(x) \le \frac{3}{4}\log^2(x+1)$ for any $x$ satisfying $\log^2(x+1) \le \frac{2}{3}$. The assertion of the lemma follows directly from these inequalities. $\square$

LEMMA 8.2. *For any sequence of positive numbers* $\theta_1, \dots, \theta_M$ *such that*

$$\sum_{i=1}^{M-1} \sqrt{\mathcal{K}(\theta_i, \theta_{i+1})} \le \frac{1}{3}$$

*it holds*

$$(8.3) \qquad \sqrt{\mathcal{K}(\theta_1, \theta_n)} \le \frac{3}{2} \sum_{i=1}^{M-1} \sqrt{\mathcal{K}(\theta_i, \theta_{i+1})}.$$

PROOF. To prove (8.3) note that by (8.1),

$$\left| \log \frac{\theta_1}{\theta_n} \right| \le \sum_{i=1}^{M-1} \left| \log \frac{\theta_i}{\theta_{i+1}} \right| \le \sqrt{3} \sum_{i=1}^{M-1} \sqrt{\mathcal{K}(\theta_i, \theta_{i+1})} \le \frac{1}{\sqrt{3}}.$$

Then using (8.2),

$$\sqrt{\mathcal{K}(\theta_1, \theta_n)} \le \frac{\sqrt{3}}{2} \left| \log \frac{\theta_1}{\theta_n} \right|,$$

which in conjunction with (8.1) proves (8.3). $\square$

LEMMA 8.3. *If the sequence* $\tau_n \ge x_0$, $n = 1, 2, \dots$, *is such that* $n_{\tau_n} \to \infty$ *as* $n \to \infty$, *then* $\widehat{n}_{\tau_n} \overset{\mathbf{P}_F}{\asymp} n_{\tau_n}$ *as* $n \to \infty$.

PROOF. By Chebyshev's exponential inequality, for any $u > 0$ and $\varepsilon \in (0, 1)$,

$$\mathbf{P}_F(\widehat{n}_{\tau_n}/n_{\tau_n} < 1 - \varepsilon) \le \exp(u(1-\varepsilon)n_{\tau_n} + n_{\tau_n}(e^{-u} - 1))$$
$$\le \exp(-u\varepsilon n_{\tau_n} + u^2 n_{\tau_n}).$$

In the same way $\mathbf{P}_F(\widehat{n}_{\tau_n}/n_{\tau_n} > 1 + \varepsilon) \le \exp(-u\varepsilon n_{\tau_n} + u^2 n_{\tau_n})$. Choosing $u = \varepsilon/2$ one gets

$$\mathbf{P}_F\left( \left| \frac{\widehat{n}_{\tau_n}}{n_{\tau_n}} - 1 \right| > \varepsilon \right) \le 2\exp\left( -\frac{\varepsilon^2}{4} n_{\tau_n} \right).$$

Since $n_{\tau_n} \to \infty$ one gets the first assertion. $\square$

LEMMA 8.4. *For any sequence* $\tau_n$, $n = 1, 2, \dots$, *satisfying* $n_{\tau_n} \to \infty$ *as* $n \to \infty$, *it holds* $\lim_{n \to \infty} \mathbf{P}_F(X_{n,k} > \tau_n) = 1$, *for any given natural number* $k$.

PROOF. By Lemma 8.3 $\mathbf{P}_F(\widehat{n}_{\tau_n} > k) \to 1$ as $n \to \infty$. Since $\mathbf{P}_F(X_{n,k} > \tau_n) = \mathbf{P}_F(\widehat{n}_{\tau_n} > k)$ we obtain the assertion of the lemma. □

LEMMA 8.5. *Assume that $P \sim Q$. Then*

$$\chi^2(P, Q) \le \mathbf{E}_P\left(\log^2 \frac{dQ}{dP} \exp\left(\left|\log \frac{dQ}{dP}\right|\right)\right).$$

PROOF. It is easy to see that $\chi^2(P, Q) = \int g(\frac{dQ}{dP})\,dP$, where $g(x) = \frac{(x-1)^2}{x}$. Since $(x-1)^2 \le e^{2\log x}\log^2 x$, for $x > 1$ and $(x-1)^2 \le \log^2 x$, for $x \in (0,1)$, we get $g(x) \le \log^2 x \exp(|\log x|)$, for $x > 0$. □

PROPOSITION 8.6. *Assume that d.f.'s $F$ and $G$ are such that it holds $\rho_t = \sup_{x \ge t} \rho_*(\alpha_F(x), \alpha_G(x)) \le \varepsilon_0$ and $\int_1^\infty (1 + \log x)^2 x^{\varepsilon_0} F_t(dx) \le \varepsilon_1$. Then $\chi^2(F_t, G_t) \le C(\varepsilon_0, \varepsilon_1)\rho_t^2$, where $C(\varepsilon_0, \varepsilon_1) = \varepsilon_1 e^{\varepsilon_0}$.*

PROOF. Since

$$\log \frac{dF_t(x)}{dG_t(x)} = \log \frac{\alpha_G(xt)}{\alpha_F(xt)} + \int_t^{xt}\left(\frac{1}{\alpha_G(u)} - \frac{1}{\alpha_F(u)}\right)\frac{du}{u}, \qquad x \ge 1,$$

it holds $|\log \frac{dF_t(x)}{dG_t(x)}| \le \rho_t(1 + \log x)$. Using Lemma 8.5, with $P = G_t$ and $Q = G_t$ one gets $\chi^2(F_t, G_t) \le \rho_t^2 e^{\rho_t} \int_1^\infty (1 + \log x)^2 x^{\rho_t} F_t(dx)$. This implies the assertion of the proposition. □

## REFERENCES

[1] BEIRLANT, J., GOEGEBEUR, YU., SEGERS, J. and TEUGELS, J. (2004). *Statistics of Extremes. Theory and Applications.* Wiley, Chichester. MR2108013

[2] BINGHAM, N. H., GOLDIE, C. M. and TEUGELS, J. L. (1987). *Regular Variation.* Cambridge Univ. Press. MR0898871

[3] DANIELSSON, J., DE HAAN, L., PENG, L. and DE VRIES, C. G. (2001). Using a bootstrap method to choose the sample fraction in tail index estimation. *J. Multivariate Anal.* **76** 226–248. MR1821820

[4] DONOHO, D. and JONSTONE, J. (1994). Ideal spatial adaptation by wavelet shrinkage. *Biometrica* **81** 425–455. MR1311089

[5] DREES, H. and KAUFMANN, E. (1998). Selecting the optimal sample fraction in univariate extreme value estimation. *Stochastic Process. Appl.* **75** 149–172. MR1632189

[6] DREES, H. (1998). Optimal rates of convergence for estimates of the extreme value index. *Ann. Statist.* **26** 434–448. MR1608148

[7] EMBRECHTS, P., KLÜPPELBERG, K. and MIKOSCH, T. (1997). *Modelling Extremal Events.* Springer, Berlin. MR1458613

[8] GOMES, M. I. and OLIVEIRA, O. (2001). The bootstrap methodology in statistics of extremes—choice of the optimal sample fraction. *Extremes* **4** 331–358. MR1924234

[9] GRAMA, I. and SPOKOINY, V. (2007). Pareto approximation of the tail by local exponential modeling. *Bul. Acad. Ştiinţe Republ. Mold. Mat.* **1** 3–24. MR2320489
[10] HALL, P. (1982). On some simple estimates of an exponent of regular variation. *J. Roy. Statist. Soc. Ser. B* **44** 37–42. MR0655370
[11] HALL, P. (1990). Using the bootstrap to estimate mean squared error and select smoothing parameter in nonparametric problems. *J. Multivariate Anal.* **32** 177–203. MR1046764
[12] HALL, P. and WELSH, A. H. (1984). Best attainable rates of convergence for estimates of parameters of regular variation. *Ann. Statist.* **12** 1079–1084. MR0751294
[13] HALL, P. and WELSH, A. H. (1985). Adaptive estimates of regular variation. *Ann. Statist.* **13** 331–341. MR0773171
[14] HILL, B. M. (1975). A simple general approach to inference about the tail of a distribution. *Ann. Statist.* **3** 1163–1174. MR0378204
[15] HUISMAN, R., KOEDIJK, K. G., KOOL, C. J. M. and PALM, F. (2001). Tail-index estimates in small samples. *J. Bus. Econom. Statist.* **19** 208–216. MR1939710
[16] LEPSKI, O. V. (1990). One problem of adaptive estimation in Gaussian white noise. *Teor. Veroyatnost. i Primenen.* **35** 459–470. MR1091202
[17] LEPSKI, O. V. and SPOKOINY, V. G. (1995). Local adaptation to inhomogeneous smoothness: Resolution level. *Math. Methods Statist.* **4** 239–258. MR1355247
[18] MASON, D. (1982). Laws of large numbers for sums of extreme values. *Ann. Probab.* **10** 754–764. MR0659544
[19] REISS, R.-D. (1989). *Approximate Distributions of Order Statistics: With Applications to Nonparametric Statistics.* Springer, New York. MR0988164
[20] RESNICK, S. I. (1997). Heavy tail modelling and teletraffic data. *Ann. Statist.* **25** 1805–1869. MR1474072
[21] WEISSMAN, I. (1978). Estimation of parameters and large quantiles based on the *k* largest observations. *J. Amer. Statist. Assoc.* **73** 812–815. MR0521329

UNIVERSITÉ DE BRETAGNE SUD
CERYC
56000 VANNES
FRANCE
E-MAIL: ion.grama@univ-ubs.fr

WEIERSTRASS-INSTITUTE
MOHRENSTR. 39
10117 BERLIN
GERMANY
E-MAIL: spokoiny@wias-berlin.de